\magnification=\magstep1
\input amstex
\UseAMSsymbols
\input pictex
\NoBlackBoxes
\pageno=1

   \font\rmk=cmr8    \font\itk=cmti8  \font\ttk=cmtt8


\def\Tr{\operatorname{Tr}}

\def\mo{\operatorname{mod}}

\def\Hom{\operatorname{Hom}}

\def\Ext{\operatorname{Ext}}

\def\ql{\operatorname{ql}}
\def\bdim{\operatorname{\bold{dim}}}
\def\minus{\!-\!}

\def\arr#1#2{\arrow <1.5mm> [0.25,0.75] from #1 to #2}

\vglue2truecm
\centerline{\bf The shift orbits of the graded Kronecker modules.}
     \bigskip
\centerline{Claus Michael Ringel}
           \bigskip\bigskip
\plainfootnote{} {\rmk 2010 \itk Mathematics Subject Classification. \rmk Primary:
		16G20, 
                16G70, 
		15A22. 
 Secondary:     05C20, 
               	16G60. 
\itk Keywords: \rmk
  Kronecker module. Regular tree. Matrix pencil.  
}
{\narrower \narrower \noindent {\bf Abstract.} Let $k$ be a field. The
Kronecker modules (or matrix pencils) are the representations of the $n$-Kronecker quiver $K(n)$;
this is the quiver with two vertices, namely a sink and a source, and $n$ arrows. The 
representations of $K(n)$ play an important role in many parts of mathematics. For $n = 2$, 
the indecomposable representations have been classified by Kronecker, but not much is known in
case $n\ge 3$. In this paper, we usually will assume that $n\ge 3$. 

The universal cover of $K(n)$ is the $n$-regular tree with bipartite orientation. 
Let $T(n)$ be the $n$-regular tree. We fix a bipartite orientation $\Omega$ of $T(n)$; the
opposite orientation will be denoted by $\sigma \Omega$ (thus $\sigma^2\Omega = \Omega$).
The $k$-representations of the quiver $(T(n),\Omega)$ can be considered as graded Kronecker
modules and we denote by $\mo (T(n),\Omega)$ the category of these graded Kronecker modules.
Only few Kronecker modules can be graded, but the graded Kronecker modules provide hints about the
behavior of general Kronecker modules. 

There is a reflection functor $\sigma\:\mo (T(n),\Omega) \to \mo (T(n),\sigma\Omega)$
(the simultaneous Bernstein-Gelfand-Ponomarev
reflection at all sinks); it will be called the {\it shift functor.}
An indecomposable graded Kronecker module $M$ is said to be {\it regular}
provided $\sigma^tM \neq 0$ for all $t\in \Bbb Z$.

If $p,q$ are vertices of $T(n)$, we denote by
$d(p,q)$ their distance. Now, let $M$ be an indecomposable regular
representation of $(T(n),\Omega)$. We attach to $M$ a positive
integer $r_0(M)$ and a pair $p(M),q(M)$ of vertices of $T(n)$ with
$0 \le d(p(M),q(M)) \le r_0(M)$ and such that
$p(M)$ is a sink if and only if $r_0(M)$ is even.
Here are the essential properties of the invariants $r_0(M),$ $p(M),q(M)$.
The $\sigma$-orbit of $M$ 
contains a unique sink module $M_0$ with smallest possible radius, say with
radius $r_0 = r_0(M)$. For $i\in \Bbb Z$, we write $M_i = \sigma^iM_0$ and call
$i = \iota(M_i)$ the {\it index} of $M_i.$ 
By duality, the
$\sigma$-orbit of $M$ contains a unique source module with radius $r_0$, say $M_{b+1}$,
and we have $b\ge 0$.
Let $p(M)$ be the center of $M_0$, let $q(M)$ be the center of $M_{b+1}$, and
denote by $(p = a_0,a_1,\dots,a_{b-1},a_b = q)$ the unique path from $p$ to $q$.
For $i\ge 0$, the module $M_{-i}$ is a sink module
with center $p(M)$ and radius $r_0+i$, whereas the module $M_{b+1+i}$ is a source module
with center $q(M)$ and radius $r_0+i$. The remaining modules $M_i$ (with $1 \le i \le b$)
are flow modules with radius $r_0-1$, and center $\{a_{i-1},a_i\}.$

By construction, the triple $r_0(M),\ p(M),\ q(M)$ is invariant under the shift.
We show that any triple $r_0,p,q$ consisting of a positive
integer $r_0$, and vertices $p,q$ of $T(n)$ with
$0 \le d(p,q) \le r_0$ and such that $p$ a sink if and only if $r_0$ is even, arises in this way.

If $M,M'$ are regular indecomposable modules with an irreducible map $M\to M'$, then we
show that $\iota(M') = \iota(M) -1$. In this way, we obtain a global way to 
index the regular indecomposable modules.
\par}
	\bigskip\bigskip
{\bf 1. Introduction.}
	     \bigskip

Let $k$ be a field. 
We denote by $T(n)$ the $n$-regular tree and usually, we will assume that $n\ge 3$. 
If $(a_0,\dots, a_t)$ is a path in $T(n)$
and $t = 2r$ is even, then $a_r$ will be called its {\it center}
and $r$ its {\it radius.}
If $t = 2r+1$, then the pair $\{a_r,a_{r+1}\}$ is called its {\it center} and $r$ its {\it radius.}

We fix a bipartite orientation $\Omega$ of $T(n)$ 
and denote by $\sigma \Omega$
the opposite orientation and we set $\sigma^2\Omega = \Omega$. (Recall that
an orientation of a graph is bipartite provided any vertex is a sink or a source;
there are just two bipartite orientations of $T(n)$.)
The finite-dimensional $k$-representations of the quivers $(T(n),\Omega)$ and $(T(n),\sigma\Omega)$ 
will be called the {\it graded Kronecker
modules} or just {\it modules} and we denote by $\mo (T(n),\Omega)$ the category of the
$k$-representations of $(T(n),\Omega)$.

Let $M$ be an indecomposable module.
We denote by $T(M)$ its support, it is by definition the full subgraph of
$T(n)$ given by all vertices $a$ with $M_a \neq 0.$
Any path in $T(M)$ of maximal length will be called a {\it diameter path}.
The length of the diameter paths will be denoted by $d(M)$ and called the
{\it diameter} of $M$.
By definition, the {\it radius} $r(M)$ of $M$ is the radius of a diameter path of $M$,
thus $r(M) = \lfloor \frac12 d(M) \rfloor$.
{\it All diameter paths of $M$ have the same center} (see section 2), called the {\it center of $M$}.

The module $M$ is called a {\it sink module} (or a {\it source module})
provided a diameter path (and hence all diameter paths)
starts and ends in sinks (or in sources, respectively). Of course, the diameter of
a sink module or a source module is even. If the diameter of $M$ is odd (so that any
diameter path for $M$ connects a sink with a source), then $M$
will be called a {\it flow module}.
The center of a sink or a source module is a vertex, the center of a flow module
is an edge.

There is a reflection functor $\sigma\:\mo (T(n),\Omega) \to \mo (T(n),\sigma\Omega)$
(the simultaneous Bernstein-Gelfand-Ponomarev
reflection at all sinks); it will be called the {\it shift functor.}
An indecomposable module $M$ is said to be {\it regular}
provided $\sigma^tM \neq 0$ for all $t\in \Bbb Z$. If $M$ is regular,
we attach to $M$ a positive integer $r_0(M)$ and a path $\pi(M) = (a_0,\dots,a_b)$ in $T(n)$,
such that $a_0$ is a sink if $r_0(M)$ is even and a source otherwise.
     \medskip
{\bf Theorem 1.} {\it Let $M$ be a regular indecomposable module.
Then the shift orbit of $M$ contains a unique sink module $M_0$
with smallest possible radius, say with radius $r_0 = r_0(M)$.
Let $M_i = \sigma^i M_0$ for all $i\in \Bbb Z$.
Then there is an integer $0\le b \le r_0(M)$ and a path
$(a_0,\dots,a_b)$ in $T(n)$ with the following properties:
      \smallskip
\item{\rm (1)}  For $i\ge 0$, the module $M_{-i}$ is a sink module with radius $r_0+i$ and
        center $a_0$.
\item{\rm (2)} For $1\le i \le b$, the module $M_i$ is a flow module with radius $r_0-1$ and
   center $\{a_{i-1},a_i\}$.
\item{\rm (3)} For $i\ge 0$, the module $M_{b+1+i}$ is a source module with radius $r_0+i$ and
     center $a_b$.
    \smallskip
\noindent
The integer $r_0(M)$ is positive. If $r_0(M)$ is even, then $a_0$ is a sink, otherwise a source.}
 	\bigskip
We call $\pi(M) = (a_0,\dots,a_b)$ the {\it center path} of $M$. 
and define $p(M) = a_0,\ q(M) = a_b$ and $b(M) = b = d(p(M),q(M)).$ By
construction, $r_0(M)$ and $\pi(M)$, thus also $p(M),q(M)$ and $b(M)$
are invariant under the shift.
    \bigskip
{\bf Theorem 2.} {\it Let $r$ be a positive integer. Let
$(a_0,\dots,a_b)$ be a path in $T(n)$ of length $0\le b \le r$.
If $r$ is even, we assume that $a_0$ is a sink, otherwise that $a_0$ is a source. Then
there is a regular indecomposable module $M$ with $r_0(M) = r$ and
$\pi(M) = (a_0,\dots,a_b).$}
    \bigskip
If $M$ is a regular indecomposable module and $M_0 = \sigma^{-i}M$ is the sink module  in the
$\sigma$-orbit of $M$ with smallest possible radius, we call $i = \iota(M)$ the {\it index} of
$M$, thus $M = \sigma^{\iota(M)}M_0.$
	\bigskip
{\bf Theorem 3.} {\it If $M, M'$ are regular indecomposable modules and $M \to M'$
is an irreducible map, then $\iota(M') = \iota(M)-1.$}
	\bigskip 
Since the regular Auslander-Reiten components of the 
category $\mo (T(n),\Omega)$ are of the form $\Bbb Z\Bbb A_\infty$, we may define
an operator $\eta$ on the set of isomorphism classes of indecomposable regular modules as follows:
if $0\to X \to Y\oplus Y' \to Z \to 0$ is an Auslander-Reiten
sequence with $Y,Y'$ indecomposable, and $|Y| < |Y'|$, let $\eta Y = Y'$.
By Theorem 3, $Y$ and $\eta(Y)$ have the same index.
	\medskip 
{\bf Theorem 4.}
{\it Let $Y$ be an indecomposable regular module. Then $Y$ is a sink module
 (or a flow module, or a source module), if and only if $\eta Y$ is a sink module
 (or a flow module, or a source module, respectively). Also, $\eta Y$ has the same center as $Y$
  and  $r(\eta Y) = r(Y)+2.$}
	\bigskip 
The proof of Theorem 1 will be given in section 3, the proof of Theorem 2 in section 5.
Section 2 collects relevant definitions and some preliminary results.
In section 4, we discuss in which way the diameter paths of $M$ and $\sigma M$ are related.

According to Gabriel,
the Auslander-Reiten translation $\tau$ in $\mo(T(n),\Omega)$ is nothing else
than $\sigma^2$ (see [Gab]). Thus, Theorem 1 provides information about 
the $\tau$-orbits of the regular
indecomposable modules and on the shape of the components of the Auslander-Reiten quiver,
see sections 6 and 7. In particular, section 7 will contain a proof of Theorem 3 and Theorem 4. 

It will be convenient to call a module {\it even} or 
{\it odd,} provided it is indecomposable and regular
and its index is even or odd, respectively.
	\medskip
{\bf Acknowledgment.} The essential parts of the paper were presented in lectures at Shanghai
and Beijing in 2015. The author is grateful to comments by mathematicians 
who attended these lectures. The author also wants to thank Daniel Bissinger, Philip Fahr, 
Rolf Farnsteiner and
Otto Kerner for many fruitful discussions concerning the behavior of 
Kronecker and graded Kronecker modules. 
    \bigskip\bigskip
\vfill\eject
{\bf 2. Preliminaries.}
   	 \bigskip
{\bf (2.1) Trees.} A {\it graph} $G$ is given by a set $G_0$ (called the
{\it vertices} of $G$) and a set
$G_1$ of subsets of $G_0$ which have cardinality $2$ (the elements of $G_1$ are called the
{\it edges} of $G$; if $\{a,a'\}$ is an edge, then $a,a'$ are also called {\it neighbors}).
A {\it path} in $G$ of length $t\ge 0$ is a sequence $(a_0,\dots,a_t)$ of
vertices such that for $1\le i\le t$, the vertices $a_{i-1}, a_i$ are neighbors, and
for $1\le i < t$, we have $a_{i-1}\neq a_{i+1}$; we say that the path $(a_0,\dots,a_t)$
connects $a_0$ and $a_t$, or also that it goes from $a_0$ to $a_t$. Instead of $\pi = (a_0,\dots,a_t)$
we also will write $\pi = (q_0|\{a_0,a_1\},\dots,\{a_{t-1},a_t\}|a_t).$
If $(a_0,\dots,a_t)$ is a path of even length $t = 2r$, then $a_r$ (or also $\{a_r\}$)
is called its {\it center} and $r$ its {\it radius.} 
If $(a_0,\dots,a_t)$ is a path of odd length $t = 2r+1$, then $\{a_r,a_{r+1}\}$
is called its {\it center} and $r$ its {\it radius.} Thus, the center of a path is either
a vertex or an edge.

A graph is
{\it connected} provided for any pair $a,a'$ of vertices, there
is a path which connects $a$ and $a'$.
A graph $G$ is a {\it tree} provided $T$ is connected and the only paths which
connect a vertex with itself are the paths of length $0$.
If $T$ is a tree, for every pair $a,a'$ of vertices, there is a unique path
going from $a$ to $a'$; the length of this path is called the {\it distance}
between $a,a'$ and is denoted by $d(a,a')$. If $x$ is a vertex of $T$ and $A$ a subset of
$T$, let $d(x,A) = \min \{d(x,a)\mid a\in A\}$.

If $T$ is a finite tree, the paths of maximal possible length $d$ are called the
{\it diameter paths} and $d$ is called the {\it diameter} of $T$.
  \medskip
{\bf Lemma 1.} {\it Let $T$ be a finite tree and $\pi$ a diameter path in $T$ with center $\{c,c'\}$
  and radius $r$. If $x$ is any vertex of $T$, then $d(x,\{c,c'\}) \le r.$}
      \medskip
      Proof. Let $\pi = (a_0,\dots,a_d)$ and $A = \{a_0,\dots,a_d\}$
      Let $x$ be a vertex of $T$. Let $j = d(x,A)$. Note that there is a path
      $(x = x_0,\dots,x_j)$ with $x_j = a_i$ for some $i$ and $x_{j-1} \neq a_{i-1}$, $x_{j-1}
      \neq a_{i+1}$.

First, assume that $i\le r$.
Since $x_{j-1} \neq a_{i+1}$, we see that
$(x_0,\dots, x_{j-1}, x_j = a_i, a_{i+1},\dots a_0)$ is a path. Its length is $j+(d-i)$.
Since $d$ is the maximal length of a path in $T$, it follows that $j\le i$. Thus,
$(x_0,\dots, x_{j-1}, x_j = a_i,\dots,a_r)$ is a path of length $j+(r-i) \le r.$
Since $a_r\in\{c,c'\}$, we have $d(x,\{c,c'\}) \le r.$

In case $i > r$, we obtain in a similar way a path from $x_0$ via $a_i$ to $c' = a_{d-r}$,
and its length is again at most $r$. \hfill$\square$
    \medskip

{\bf Corollary.}
{\it All diameter paths of a finite tree $T$ have the same center} (and, of course, the same radius).
     \medskip

Proof: Let $\pi = (a_0,a_1,\dots,a_d)$ be a diameter path with center $\{c,c'\}$ and radius $r$.
Consider a diameter path $\pi'$ which connects $x$ with $x'$. Then $d(x,\{c,c'\}) \le r$ and
$d(x',\{c,c'\}) \le r$. Thus there is a path from $x$ say to $c$ of length at most $r$.
If there is also a path from $x'$ to $c$ of length at most $r$, then $d(x,x') \le 2r$, thus
$d = d(x,x')\le 2r\le d$. It follows that $d$ is even, thus $c=c'$, and that the path from $x$ to $x'$
runs through $c$. But this means that $c$ is the center of $\pi'$.

Thus, it remains to consider the case that $c\neq c'$ and that there is a path from $x'$ to $c'$
of length at most $r$. Then $\pi'$ has to be the concatenation of the path from $x$ to $c$ with
the edge $\{c,c'\}$ and the path from $c'$ to $x'$, thus the center of $\pi'$ has to be $\{c,c'\}$.
 \hfill$\square$
 \medskip

By definition, the {\it center} $C(T)$ and the {\it radius} $r(T)$
of $T$ are the center and the radius of the diameter paths of $T$.

   \medskip
If $T$ is a finite tree, the {\it boundary} of $T$ is defined to be 
the set of endpoints of the diameter paths (these are leaves of $T$; but usually, $T$ will have
additionsl leaves).
       \medskip
A vertex $x$ of a graph is called a {\it leaf} provided $x$ has at most one neighbor.
Of course, the boundary vertices of a finite tree are leaves.

       \bigskip 
{\bf (2.2) Orientations.} Given a graph $G$, an {\it orientation} $\Omega$ of $G$ is
a map $\Omega\:G_1 \to G_0\times G_0$ such that $\Omega(\{a,a'\}$ is either $(a,a')$ or
$(a',a)$.
If $\Omega(\{a,a'\}) = (a,a')$, we write $a\to a'$ and call $a$ the
{\it start} and $a'$
the {\it target} of the arrow $a\to a'$. A vertex $a$ of $(G,\Omega)$ is a {\it sink} (or a
{\it source}) provided $a$ is not the start (or the target, respectively) of any arrow.
The orientation $\Omega$ is called {\it bipartite} provided any vertex is a sink or a source.
If $\Omega$ is an orientation of $G$, then $(G,\Omega)$ is called an
{\it oriented graph.}

Let $(T,\Omega)$ be a finite tree with bipartite orientation.
If the diameter of $T$ is odd, then any diameter
path connects a sink with a source. If the diameter of $T$ is even, and there is a
diameter path which connects a sink with a sink, then all diameter paths connect sinks with
sinks (since all diameter paths have the same center).
      \bigskip
      {\bf (2.3) Balls.} Now we consider finite subgraphs of $T(n)$.
      If $c_1,\dots,c_t$ are vertices of $T(n)$, let $B_r(c_1,\dots,c_t)$
      be the full subgraph of
      $T(n)$ given by all vertices $a$ with $d(a,c_i)\le r$ for some $1\le i \le t.$
      We are only interested in the case $t=1$ and in the case $t=2$ with $c_1,c_2$ being neighbors.
      Of course, given a vertex $c$, then $B_r(c)$ has center $c$ and radius $r$.
      Similarly, if $c_1,c_2$ are neighbors, then $B_r(c_1,c_2)$ has center $\{c_1,c_2\}$
      and radius $r$. We call $B_r(c)$ and $B_r(c_1,c_2)$ the {\it ball} with radius $r$
      and with center $c$ or $\{c_1,c_2\}$, respectively.
      	  \bigskip
	  {\bf (2.4) Representations.} If $G$ is a graph and $\Omega$ an orientation of $G$, the
	  oriented graph $(G,\Omega)$ is nothing also than a quiver without loops and multiple arrows,
	  and we denote by $\mo(G,\Omega)$ the category of finite-dimensional
	  $k$-representations of this quiver. The case we are interested in is $G = T(n)$ and
	  $\Omega$ a fixed bipartite partition. The finite-dimensional $k$-representations of $(T(n),\Omega)$
	  are the {\it modules} we are dealing with.

If $M$ is an indecomposable module, we denote by $T(M)$ its support; it is the full
subgraph given by all vertices $a$ of $T(n)$ with $M_a \neq 0$. Of course, this is a finite
tree and we write $C(M) = C(T(M))$ and $r(M) = r(T(M))$ and call it the center and the
radius of $M$, respectively.

An indecomposable module will be called a
{\it sink} module provided any (and thus all)  diameter paths in $T(M)$ start and end in sinks,
and a
{\it source} module provided any (and thus all)  diameter paths in $T(M)$ start and end in sources.
An indecomposable module with even diameter is either a sink module or a source module.
An indecomposable module with odd diameter will be called a {\it flow} module.
Note that any diameter path of a flow module connects a sink with a source.
Recall that the center of an indecomposable module is either a vertex or an edge. For the sink and
the source modules, the center is a vertex, for the flow modules, the center is an edge.
    \medskip
    If $M$ is an indecomposable module with center $C$ and radius $r$, then
    we write $B(M) = B_r(C)$. We have $T(M) \subseteq B(M)$, and $T(M)$ and $B(M)$
    have the same radius and the same center. Note that the boundary vertices of $T(M)$
    are just the leaves of $B(M)$ which belong to $T(M)$.
    	\bigskip
{\bf (2.5) Completeness of an indecomposable module.}
Let $M$ be a sink or a source module with center $C$ and radius $r$.
The module $M$ is said to be {\it complete} provided $r \ge 1$ and
for any path $(x(0), x(1),\dots, x(r))$ in $T(n)$ with $x(r) = p$ and such that
all the vertices $x(i)$ with $1\le i \le r$ belong to $T(M),$ one has
$\dim M_{x(0)} = \dim M_{x(1)}$ (so that, in particular, also $x(0)$ belongs to $T(M)$). 

Let $M$ be a flow module with radius $r$.
The module $M$ is said to be {\it complete} provided $r\ge 1$ and such that
for any path $(x(0), x(1),\dots, x(r))$ in $T(n)$ such that
$x(r)$ is a central vertex, but $x(r-1)$ is not a central vertex and such that
all the vertices $x(i)$ with $1\le i \le r$ belong to $T(M),$ one has
$\dim M_{x(0)} = \dim M_{x(1)}$ (so that, in particular, also $x(0)$ belongs to $T(M)$).
     \medskip
A module is said to be {\it incomplete} provided it is indecomposable
and not complete.
Let us stress that by definition
the simple modules and the indecomposable modules of length $2$
(these are the indecomposable modules
with radius $0$) are incomplete.
            \bigskip
{\bf (2.6) The shift functor.}
If $x$ is a sink of a quiver $\Delta$, we denote by $\sigma_x\Delta$ the quiver
obtained from $\Delta$ by changing the orientation of all the arrows of $\Delta$ ending in $x$;
in particular, $x$ is a source of $\sigma_x\Delta$. Bernstein, Gelfand and Ponomarev ([BGP]) have
defined a so-called reflection functor $\sigma_x:\mo \Delta \to \mo \sigma_x\Delta$;
We denote by $\sigma_x^-$ the left adjoint of $\sigma_x$. Note that $\sigma_x$ sends $S(x)$
to zero and if $M$ is an indecomposable representation of $\Delta$ different from $S(x)$, then
$\sigma_x^-\sigma_x(M)$ is naturally isomorphic to $M$.

The reflection functors $\sigma_x$ and $\sigma_{x'}$ for different sinks $x,x'$
commute, thus the composition $\sigma$ of the reflection functors $\sigma_x$ for all sinks $x$ of
$\Delta$ is well-defined and independent of the order used, it is a functor
$\mo\Delta$ to $\mo\sigma\Delta$, where $\sigma\Delta$ is the quiver obtained from
$\Delta$ by changing the orientation of all the arrows of $\Delta$ ending in a sink.
We denote by $\sigma^-$ the left adjoint of $\sigma$. Note that $\sigma$ sends $S(x)$ with
$x$ any sink to zero. If $M$ is an indecomposable representation of $\Delta$ which is not simple
projective, then $\sigma^-\sigma(M)$ is naturally isomorphic to $M$.

In particular, if $\Delta = (T(n),\Omega)$, then $\sigma\Delta = (T(n),\sigma\Omega)$,
where $\sigma\Omega$ is the opposite orientation (the orientation of all arrows is changed).
The reflection functor $\sigma\:\mo (T(n),\Omega) \to \mo (T(n),\sigma\Omega)$
will be called the {\it shift functor.}

Of particular interest is the composition $\sigma^2\:\mo (T(n),\Omega) \to \mo (T(n),\Omega)$.
According to Gabriel [Gab], this is just the Auslander-Reiten functor $\tau = D\Tr$ (``dual of transpose'');
note that [Gab] clarified and corrected a previous assertion by Brenner and Butler [BB].
     \bigskip
{\bf (2.7)} From now on, we consider finite-dimensional $k$-representations
of the quiver $(T(n),\Omega)$ (or also of $(T(n),\sigma\Omega)$), 
thus {\it graded Kronecker
modules} and call them just {\it modules}.

An indecomposable module $M$ will be said to be {\it preprojective} provided
$\sigma^tM = 0$ for some $t\ge 1.$ The preprojective modules are of the form
$P_{2s}(x)= \sigma^{-2s}S(x)$ for some sink $x$ in $(T(n),\Omega)$ or of the form
$P_{2s+1}(y) = \sigma^{-2s-1}S(y)$ for some sink $y$ in $(T(n),\sigma\Omega)$, with $s\ge 0$.
Note that the preprojective modules are sink modules
(these modules have been considered in detail in [FR1]).
Dually, an indecomposable module $M$ will be said to be {\it preinjective} provided
$\sigma^{-t}M = 0$ for some $t\ge 1.$ The preinjective modules are of the form
$Q_{2s}(y)= \sigma^{2s}S(y)$ for some source $y$ in $(T(n),\Omega)$ or of the form
$Q_{2s+1}(x) = \sigma^{2s+1}S(x)$ for some source $x$ in $(T(n),\sigma\Omega)$, with $s\ge 0$.
The preinjective modules are source modules.
The preprojective and the preinjective modules are simple or complete.

An indecomposable module $M$ will be called {\it regular} provided it is neither
preprojective, nor preinjective, thus provided $\sigma^tM \neq 0$ for all $t\in \Bbb Z.$
     \bigskip\bigskip

{\bf 3. Proof of Theorem 1.}
   \bigskip
   Let $M$ be a regular indecomposable module. We apply $\sigma$ and describe the change of the
   diameter. Six different cases (1) \dots (6) will occur.
   	     \medskip
{\it Let $M$ be a sink module with center $c$. Then:
$$
d(M)\minus 2 \le d(\sigma M) \le d(M).
$$
There are the following three possibilities:}
      \smallskip
{\bf (1)} {\it $d(\sigma M) = d(M)\minus 2$
if and only if $\sigma M$ is a sink module
if and only if $M$ is complete.
 In this case,
the center of $\sigma M$ is $c$.}
    \smallskip
{\bf (2)} {\it  $d(\sigma M) = d(M)\minus 1$ if and only if $\sigma M$ is a flow module.
The center of $M$ is of the form $\{c,c'\}$.
 with a path $(x_0,x_1,\dots,x_r)$, where $x_0 = c,
 x_1 = c'$ such that $x_r$ is a source for $\sigma \Omega$.}
     \smallskip
{\bf (3)} {\it  $d(\sigma M) = d(M)$ if and only if $\sigma M$ is a source module.
In this case, the center of $\sigma M$ is $c$ again.}
   \medskip

Proof.  Let $M$ be a sink module with center $c$ and radius $r$, thus
$T(M) \subseteq B_r(c)$.
Since the leaves of $B_r(c)$ are sinks, $\sigma(B_r(c)) \subseteq B_r(c).$ This shows
that $d(\sigma M) \le d(M)$.
Given a diameter path $\pi$
in $T(M)$ say connecting the leaves $x$ and $x'$, the vertices of $\pi$ different from $x$
and $x'$ are in the support of $\sigma M$. This shows that $d(\sigma M) \ge d(M)\minus 2.$
This completes the proof of the inequalities.

For the proof of (1) we note that $M$ is complete if and only if
$T(\sigma M) \subseteq B_{r-1}(c)$, thus if and
only if $d(\sigma M) = d(M)\minus 2.$ In this case $\sigma M$ is a sink module.
On the other hand, if there is a vertex $x \in T(\sigma M)$ which does not belong to
$B_{r-1}(c)$, then there is a diameter walk in $T(\sigma M)$ which starts at $x$
and this is a source for $\sigma \Omega$, thus $\sigma M$ is not a source module. This
completes the proof of (1).

(3) If $d(\sigma M) = d(M)$, then there is a diameter walk for $\sigma M$
starting and ending in boundary vertices of $B_r(c)$, thus in sources of $\sigma Omega$.
This shows that $\sigma M$ is a source module. On the other hand, if $\sigma M$ is
a source module, consider a diameter walk $\pi$ for $\sigma M$. Its ends cannot lie in
$B_{r-1}(c)$, thus its length is $2r = d$. In this way, we see (3).

The first part of assertion (2) follows from the inequalities and (1) and (3), since
any indecomposable module is a sink module, a source module, or a flow module.
It remains to describe in this case the center of $\sigma M$.
We start with a diameter path $(a_1,\dots,a_{2r})$ of $\sigma M$ and we may assume that
$a_1$ is a sink and $a_{2r}$ a source with respect to $\sigma \Omega$. Then the center
of $\sigma M$ is $\{a_r,a_{r+1}\}$. Let $a_0$ be a neighbor of $a_1$ different from $a_2$.
Then $(a_0,a_1,\dots,a_{2r-1})$ is a path in the support of $M = \sigma^{-1}\sigma M$
and $a_0$ is a sink with respect to $\Omega.$ Since the support of $M$ lies in $B_r(c)$,
it follows that $a_r = c$. Thus, we obtain the required path $(x_0,\dots,x_r)$ by
looking at the vertices $x_i = a_{i+r}$ for $0 \le i \le r.$
	\hfill$\square$
	\bigskip
{\bf Corollary.} {\it If $M$ is a regular indecomposable sink module with radius $r$.
	Then there is $1\le  i \le r$ such that $\sigma^i M$ is a flow module or a source module.}
	     \medskip
	     Proof. We use induction on $r$. We have $r\ge 1$ since a sink module with radius $r = 0$
	     is simple, thus not regular. Thus, we start with a sink module $M$ with radius $r\ge 1$.
	     Either $\sigma M$ is a flow or a source module (then we take $i=1$)
	     or else $\sigma M$ is a regular sink module
	     with radius $r-1$. In the latter case, induction shows that $\sigma^j (\sigma M)$ is a flow
	     or a source module, for some $1 \le j \le r-1.$ Thus $i = j+1 \le r$ and $\sigma^iM$ is a flow or
	     a source module.
	       \hfill$\square$
		\bigskip 
{\it Let $M$ be a flow module. We assume that
$(a_0,\dots,a_d)$ with $d = 2r+1$ is a diameter path of $M$
with $a_0$ a sink and $a_d$ a source. Then:
$$
  d(M) \le d(\sigma M) \le d(M)+1.
$$
There are the following two possibilities:}
      \smallskip
{\bf (4)} {\it $d(\sigma M) = d(M)$, thus $\sigma M$ is a flow module. Let $a_{d+1}$ be a neighbor
  of $a_d$ different from $a_{d-1}$. Then $(a_1,\dots,a_{d+1})$ is a diameter path for $\sigma M$.
    In particular, the center of $\sigma M$ is $\{a_{r+1},a_{r+2}\}$} (and we note that $a_1$ is a
      sink for $\sigma\Omega$ and $a_{d+1}$ is a source for $\sigma\Omega.$)
      \smallskip
{\bf (5)} {\it $d(\sigma M) = d(M)+1$, then $\sigma M$ is a source module and
  the center of $\sigma M$ is $a_{r+1}$.}
    \medskip
        Proof: Let $(a_0,\dots,a_d)$ be a diameter path of $M$ with $a_0$ a sink and $a_d$ a source.
	Let $a_{d+1}$ be a neighbor of $a_d$ different from $a_{d-1}$. Clearly, the path $(a_1,\dots,a_d,a_{d+1})$
	is in the support of $\sigma M$. This shows that $d(\sigma M) \ge d(M)$.
	Also, if $d(\sigma M) = d(M)$,  then $(a_1,\dots,a_d,a_{d+1})$ is a diameter path for $\sigma M$,
	thus the center of $\sigma M$ is $\{a_{r+1},a_{r+2}\}$. This completes the proof of (4).

Now assume that $d(\sigma M) > d(M)$. The support of $\sigma M$ is contained in $B_{r+1}(a_{r+1})$,
thus $d(\sigma M) \le 2(r+2) = d+1.$ It follows that $d(\sigma M) = d(M)+1$ and that $a_{r+1}$
is the center of $\sigma M$. The path
$(a_{r+1},\dots,a_d)$ in $(T(n),\Omega)$ ends in a source, thus the path
$(a_{r+1},\dots,a_{d+1})$ in $(T(n),\sigma\Omega)$ ends also in a source. This shows that the
boundary vertices of $B_{r+1}(a_{r+1})$ are sources, therefore $\sigma M$ is a source module.
This completes the proof of (5). 
     \hfill$\square$
		   \bigskip

{\bf (6)} {\it Let $M$ be a source module with center $c$ and radius $r$.
        Then $\sigma M$ is a source module with center $c$ and radius $r+1$.}
	     \medskip
Proof: The support of $M$ is contained in $B_r(c)$, thus the support of $\sigma M$
is contained in $B_{r+1}(c)$. It follows that the radius of $\sigma M$ is at most $r+1$.
Any diameter walk in $M$ starts and ends in a source. In $\sigma M$
such a walk is prolonged on both sides by an edge, thus we obtain in this way a
walk of length $2r+2.$  As a consequence, the radius of $\sigma M$ is precisely $r+1$.
The boundary vertices of $B_{r+1}(c)$ are sources for $\sigma \Omega$, thus $\sigma M$
is a source module with center $c$.
   \hfill$\square$
          \bigskip

In order to see that any $\sigma$-orbit contains only finitely many
flow modules, we need a new invariant.
A path which connects two (not necessarily different) source leaves will be called
a {\it source path.} Note that the length $|\pi|$ of a source path $\pi$
is an even natural number (as before, we call $\frac12|\pi|$ the {\it radius} of $\pi$).
     \medskip

{\bf Lemma 2.} {\it  Let $M$ be indecomposable and regular with radius $r(M)$
such that $T(M)$ contains source leaves. Let $w(M)$ be the maximal radius of
a source path of $M$.
Then $w(M) \le r(M)$ and $\sigma^{r(M) - w(M) +1}M$ is a source module.}
     \medskip
Proof. If $M$ is a sink or a source module, then $2r(M)$ is the diameter of $M$,
thus the maximal length of the paths in $T(M)$ and therefore $2w(M) \le 2r(M)$, thus $w(M) \le r(M)$.
If  $M$ is a flow module with diameter $d$, then $2r(M) = d-1$ is the maximal length of a path
which connects two sources, thus again $2w(M) \le 2r(M)$.

Now assume that $M$ is indecomposable and regular, with a source leaf $x$ in its support.
Then $\sigma M$ is again indecomposable and regular, with a source leaf in its support
(namely the neighbors of $x$ which are not vertices of $T(M)$ are source leaves in $T(\sigma M)$).
By induction, all the
modules $\sigma^i M$ with $i\ge 0$ are indecomposable and regular with source leaves
in the support.

Now $T(M)$ has a source path of maximal length $2w(M)$, thus $T(\sigma M)$ has a source
path of length $2w(M)+2$. Therefore $w(\sigma M) \ge w(M)+1.$ Using induction, we see that
$$
 w(\sigma^i M) \ge w(M)+i
$$
for all $i \ge 0.$

Let $t = r(M)-w(M)+1$ and assume that $\sigma^t M$ is not a source module. Then, according to (6),
none of the modules $M,\sigma M, \dots, \sigma^t M$
is a source module. Thus all these modules are sink or flow modules, and therefore
$$
 r(M) \ge r(\sigma M) \ge r(\sigma^2 M) \ge \cdots \ge r(\sigma^t M),
$$
according to (1), (2) and (4).

Altogether, we see that
$$
\align
  r(M) \ge &\ r(\sigma^{r(M)-w(M)+1}M) \cr
    \ge &\    w(\sigma^{r(M)-w(M)+1} M) \cr
      \ge  &\ w(M)+ (r(M)-w(M)+1) = r(M) + 1,
\endalign
$$
a contradiction.
\hfill$\square$
	\bigskip
{\bf Corollary 1.} {\it The shift orbit of a regular indecomposable module
contains source modules and sink modules.}
	\medskip
Proof. Consider the shift orbit $\Cal O$ of the regular indecomposable module $M$.

First, assume that $\Cal O$ contains a flow module $M'$. Then $T(M')$ has
a source leaf, thus Lemma 2 asserts that the $\Cal O$ contains a
source module. By duality, $\Cal O$ also contains a sink module.

Second, assume that the $\Cal O$ contains a sink module $M'$. Corollary 1 asserts that
one of the modules $\sigma^i M'$ with $i > 0$ has to be a flow module or a source module.
But we know already that the existence of a flow module in $\Cal O$ implies that
there is also a source module in $\Cal O$. Thus we see that $\Cal O$
contains a source module.

By duality, the existence of a sink module in $\Cal O$ implies the existence
of a source module in $\Cal O$.
   \hfill$\square$
	\bigskip
We see: Let $\Cal O$ be the shift orbit of a regular indecomposable module.
Let $M$ be a sink module in $\Cal O$.
According to corollary 1, there is some $i\ge 0$ such that $\sigma^{i+1} M$ is a flow or a source
module. Choose $i$ minimal with this property and let $M_0 = \sigma^i M.$ Then $M_0$
is an incomplete sink module. We let $M_i = \sigma^i M_0$ for all $i\in \Bbb Z$.

Then all the module $M_i$ with $i \le 0$ are sink modules (by the dual assertion of (6)).

Let $b \ge 0$ be minimal such that $M_{b+1}$ is a source module. Then
the modules $M_{b+i}$ with $i \ge 1$ are source module (by (6)). On the other hand,
all the modules $M_i$ with $1\le i \le b$ have to be flow modules.
    \bigskip

{\bf Corollary 2.} {\it If $M$ is a module with radius $r$, 
then the $\sigma$-orbit of $M$ contains at most  $r$ flow modules.}
     	        \medskip
Proof.
We can assume that the $\sigma$-orbit $\Cal O$ of $M$ contains at least one flow module $M$.
According to corollary 1, the $\sigma$-orbit of $M$ contains source and sink module.
Thus the orbit $\Cal O$ contains a
source module $M_0$ such that $M_1 = \sigma M_0$ is a flow module. We apply Lemma 2 to $M_1$
and see that there is some $1\le b \le r(M_1)$ such that the modules $M_1,\dots,M_b$ are
all the flow modules in the $\sigma$-orbit of $M$.

The assertions (1) to (6) show that $r(M_1) \le r(M)$. According to the 
Lemma 2, there is some $0 \le i \le r(M_1)$ such that
$\sigma^iM_1$ is a flow module, but $\sigma^{i+1}M_1$ is a source module.
	     \hfill$\square$
     \bigskip\bigskip
Proof of Theorem 1.
Consider a $\sigma$-orbit $\Cal O$ of regular indecomposable modules. As we
have seen, $\Cal O$ contains both a sink module $M$ and a source module, say $\sigma^s M$
for some $s\in \Bbb Z$. According to (6), $s$ must be a positive integer. Thus there is
some $0\le j < s$ such that $M_0 = \sigma^j M$ is a sink module, whereas $\sigma M_0$
is not a sink module. Thus $M_0$ is incomplete. Let $p$ be the center of $M_0$ and let $r$
be its radius.

We write $M_i = \sigma^i M_0$ for all $i\in \Bbb Z$. According to the dual assertion of (6),
the modules $M_i$ with $i\le 0$ all are sink modules with center $p$ and radius $r+i.$
According to (1), (2) and (4), the modules $M_i$ with $i > 0$ are
flow or source modules.
It follows that $M_0$ is the only incomplete sink module in $\Cal O$.

Since the number of flow modules in $\Cal O$ is finite, there is a smallest number $b\ge 0$
such that $M_{b+1}$ is a source module. Then the modules $M_i$ with $1\le i \le b$ are flow
modules. According to (6), the modules $M_i$ with $i\ge b+1$ are source modules.
This shows that $b \ge 0 $ is the number of flow modules in $\Cal O$. 

According to (2) and (4), the modules $M_i$ with $1\le i \le b$ (the flow modules in $\Cal O$)
all have radius $r-1$. According to (3), (5) and (6), the modules $M_{b+1+i}$ with $i\ge 0$
(the source modules in $\Cal O$) have radius $r+i.$ Also, according to (6), all these
source modules have the same center, say $c_b$. In case $b = 0$, we know from (3) that $c_0 = c_b$.
In case $b > 0$, we use (2), (4) and (5) in order to see that the centers of the modules
$M_0, M_1,\dots, M_{b+1}$ are of the form $c_0, \{c_0,c_1\}, \dots, \{c_{i-1},c_i\},\dots,
\{c_{b-1},c_b\}, c_b$, where $(c_0,c_1,\dots,c_b)$ is a path.
\hfill$\square$

Thus, any regular $\sigma$-orbit looks as follows:
$$
\hbox{\beginpicture
  \setcoordinatesystem units <1.05cm,.7cm>
\put{$\cdots$} at -0.7 0
\put{$M_{-2}$} at 0 0
\put{$M_{-1}$} at 1 0
\put{$M_{0}$} at 2 0
\put{$M_{1}$} at 3 0
\put{$\cdots$} at 4 0
\put{$M_{i}$} at 5 0
\put{$\cdots$} at 6 0
\put{$M_{b}$} at 7 0
\put{$M_{b+1}$} at 8.05 0
\put{$M_{b+2}$} at 9 0
\put{$M_{b+3}$} at 10 0
\put{$\cdots$} at 10.7 0
\setdashes <1mm>
\plot 2.5 1.2 2.5 -2.4 /
\plot 7.56 1.2 7.56 -2.4 /
\put{sink modules} [r] at 1.7 1
\put{flow modules} at 5 1
\put{source modules} [l] at 8.2 1

\put{module:} at -1.7 0
\put{radius:} at -1.7 -1
\put{center:} at -1.7 -2
\put{$\cdots$} at -.7 -1
\put{$\ssize r\!+\!2$} at 0 -1
\put{$\ssize r\!+\!1$} at 1 -1
\put{$\ssize r$} at 2 -1
\put{$\ssize r\minus 1$} at 3 -1
\put{$\cdots$} at 4 -1
\put{$\ssize r\minus 1$} at 5 -1
\put{$\cdots$} at 6 -1
\put{$\ssize r\minus 1$} at 7 -1
\put{$\ssize r$} at 8 -1
\put{$\ssize r\!+\!1$} at 9 -1
\put{$\ssize r\!+\!2$} at 10 -1
\put{$\cdots$} at 10.7 -1

\put{$\cdots$} at -.7 -2
\put{$c_0$\strut} at 0 -2
\put{$c_0$\strut} at 1 -2
\put{$c_0$\strut} at 2 -2
\put{$\ssize\{c_0,c_1\}$\strut} at 3 -2
\put{$\cdots$} at 4 -2
\put{$\ssize\{c_{i-1},c_i\}$\strut} at 5 -2
\put{$\cdots$} at 6 -2
\put{$\ssize\{c_{b-1},c_b\}$\strut} at 6.95 -2
\put{$c_b$\strut} at 8 -2
\put{$c_b$\strut} at 9 -2
\put{$c_b$\strut} at 10 -2
\put{$\cdots$} at 10.7 -2
\endpicture}
$$
Let us stress that there may be no flow modules in $\Cal O$, this is the case $b = 0$.
In this case, $c_0 = c_b$ and the center path is a path of length $0$.
	\bigskip
Given an indecomposable regular module $M$, we define its {\it index} $\iota(M)$
as the integer $t$ such that $\sigma^{-t}M$ is an incomplete sink module. Thus
$\iota(M) = t$ means that $M = \sigma^t M_0$ for some incomplete sink module (namely for
$M_0 = \sigma^{-t}M$).
	\medskip 
There is the following Corollary to Theorem 1.
	\medskip 
{\bf Corollary.} {\it Let $M$ be a sink module with center $p$, such that $\sigma^{b+1}M$
is a source module with center $q$. If $d(p,q) = b$, then $\iota(M) = 0.$}
   \bigskip\bigskip
{\bf 4. Diameter paths and boundary vertices.}
     \medskip
We are going to look in which way diameter paths of $M$ and of $\sigma M$ are related.
Given an indecomposable module $M$, we denote by $\gamma(M)$ the number of diameter paths for $M$.

The diameter paths may be quite different in case both $M$ and $\sigma M$ are incomplete.
Otherwise, the diameter paths of $M$ and of $\sigma M$ are very similar --- this concerns
the cases where $M$ and $\sigma M$ both are sink modules, or both are flow modules,
or both are source modules, thus, the cases (1), (4) and (6).
   \bigskip

{\bf Proposition 1.} {\it Let $M$ be a representation of $(T(n),\Omega).$}

(a) {\it Assume that $M$ and $\sigma M$ both are sink modules. Let $(a_0,\dots,a_d)$ be a path in
$T(n)$. Then $(a_0,\dots,a_d)$ is a diameter path for $M$ if and only if $(a_1,\dots,a_{d-1})$ 
is a diameter path for $\sigma M$. Thus}
$$
 \gamma(M) = (n-1)^2\gamma(\sigma M).
 $$

(b) {\it Assume that $M$ and $\sigma M$ both are flow modules. Let $(a_0,\dots,a_d)$ be a path in
$T(n)$ with $a_0$ (and also $a_d$) being sinks for $\Omega$. 
Then $(a_0,\dots,a_{d-1})$ is a diameter path for $M$ if and only if $(a_1,\dots,a_{d})$ 
is a diameter path for $\sigma M$.
Thus}
$$
 \gamma(M) = \gamma(\sigma M).
 $$

(c) {\it Assume that $M$ and $\sigma M$ both are source modules. Let $(a_0,\dots,a_{d})$ be a path in
$T(n)$.
Then $(a_1,\dots,a_{d-1})$ is a diameter path for $M$ if and only if $(a_0,\dots,a_{d})$
is a diameter path for $M$.
Thus}
$$
 \gamma(\sigma M) = (n-1)^2\gamma(M).
 $$
	\bigskip
	Proof. This follows directly from the considerations in section 3.
	\hfill$\square$
		\bigskip

{\bf Example} {\it of a sink module $M$ and a source module $M' = \sigma M$ with arbitrarily large
radius $r$ such that no edge of $T(n)$ belongs both to a diameter of $M$ and a diameter of $M'$.}

We need $n \ge 4$. For $r = 1$, let $M$ be any 3-dimensional sink module. Then $ \sigma M$ is
an $(n-1)$-dimensional source module.

Now assume that $r\ge 2.$
Let $a_0$ be a sink of $(T(n),\Omega)$ and $(a_0,a_1,\dots,a_d)$ a path of
length $d = 2r.$ Since $n\ge 4$, there is a path $(b_0,\dots,b_d)$ such that $b_r = a_r$ is the
only common vertex of these paths. Let $M$ be the thin indecomposable module with support
the full subquiver with vertices $a_j,b_j$, where $1\le j \le d-1$ and all the neighbors of
$a_1$ and $a_{d-1}$. Then $M$ is a sink module with center
$c = a_r = b_r$. The diameter
paths of $M$ are the paths of the form
$(a'_0,a'_1,\dots,a'_{d-1},a'_d)$, where $a'_j = a_j$ for $1\le j \le d-1$ (the number of
diameter paths is $(n-1)^2).$
The diameter paths of $M' = \sigma M$ are all the
paths of the form $(b'_0,\dots,b'_d)$ with $b'_j = b_j$ for $1\le j \le d-1$ (the number of
diameter paths of $M'$ is again $(n-1)^2,$ but there is no natural correspondence between
the diameter paths of $M$ and of $M'$).

For example, for $n=4$ and $r = 2$, the support of $M$ and $M'$ looks as follows:
$$
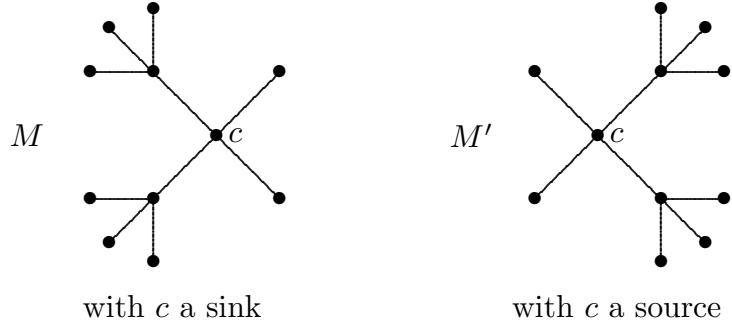
\hbox{\beginpicture
  \setcoordinatesystem units <.7cm,.7cm>
\put{\beginpicture
\put{$M$}  at -3 0
\multiput{$\bullet$} at 0 0  -1 1  -1 -1  1 1  1 -1
    -2 1  -1.7 1.7  -1 2  -2 -1  -1.7 -1.7  -1 -2 /
\plot -1.7 -1.7  1 1 /
\plot -1.7 1.7  1 -1 /
\plot -1 -2  -1 -1  -2 -1 /
\plot -2 1  -1 1  -1 2 /
\put{$c$} at 0.3 0
\put{with $c$ a sink} at -.7 -2.7
\endpicture} at 0 0

\put{\beginpicture
\put{$M'$}  at -2 0
\multiput{$\bullet$} at 0 0  -1 1  -1 -1  1 1  1 -1
        2 1  1.7 -1.7  1 2  2 -1  1.7 1.7  1 -2 /
\plot 1.7 1.7  -1 -1 /
\plot 1.7 -1.7  -1 1 /
\plot 1 -2  1 -1  2 -1 /
\plot 2 1  1 1  1 2 /
\put{$c$} at 0.3 0
\put{with $c$ a source} at .3 -2.7
\endpicture} at 7 0
\endpicture}
$$
	\bigskip
Let us add some comments concerning the number of boundary vertices of sink 
modules.
If $M$ be an indecomposable module, let $\beta(M)$ be the number of boundary vertices of $T(M)$.
Note that any ball of radius $r$ has precisely $\beta_r = n(n-1)^{r-1}$ boundary vertices.
Thus, if $M$ is a sink module (or a source) module with radius $r > 0$, then
$$
 2 \le \beta(M) \le \beta_r.
$$

{\it If $M$ is a source module, then for all $t\ge 0$}
$$
 \beta(\sigma^t M) = (n-1)^t\beta(M).
$$
This shows that the number of boundary vertices growths exponentially when we apply $\sigma$.
	\bigskip\bigskip 

{\bf 5. Examples, in particular proof of Theorem 2.}
     	\medskip
Let $(p|\gamma_1,\dots,\gamma_r|q)$ be a path,
say with arrows $\gamma_i$
between $a(i-1)$ and $a(i)$. We assume that the path starts at the sink $p = a(0)$
in case $r$ is even (see the left picture), and at the source $p = a(0)$ in case
$r$ is odd (see the right picture):
$$
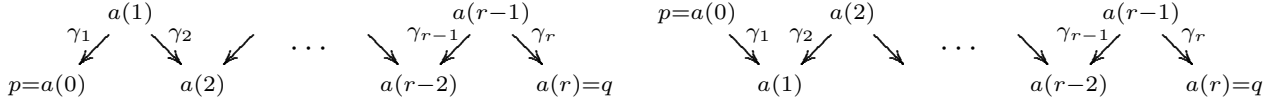
\hbox{\beginpicture
  \setcoordinatesystem units <.8cm,.8cm>
  \put{\beginpicture
  \put{$\ssize p=a(0)$} at -.15 0
\put{$\ssize a(1)$} at 1 1
\put{$\ssize a(2)$} at 2 0
\put{$\cdots$} at 3.5 0.5
\put{$\ssize a(r-2)$} at 5 0  
\put{$\ssize a(r-1)$} at 6 1
\put{$\ssize a(r)=q$} at 7.15 0
\arr{0.7 0.7}{0.3 0.3}
\arr{1.3 0.7}{1.7 0.3}
\arr{2.7 0.7}{2.3 0.3}

\arr{4.3 0.7}{4.7 0.3}
\arr{5.7 0.7}{5.3 0.3}
\arr{6.3 0.7}{6.7 0.3}
\put{$\ssize \gamma_1$} at 0.3 0.7
\put{$\ssize \gamma_2$} at 1.7 0.7
\put{$\ssize \gamma_{r-1}$} at 5.2 0.7
\put{$\ssize \gamma_r$} at 6.73 0.7
\endpicture} at 0 0
\put{\beginpicture
\put{$\ssize p = a(0)$} at -.15 1
\put{$\ssize a(1)$} at 1 0
\put{$\ssize a(2)$} at 2 1
\put{$\cdots$} at 3.5 0.5
\put{$\ssize a(r-2)$} at 5 0  
\put{$\ssize a(r-1)$} at 6 1
\put{$\ssize a(r)=q$} at 7.15 0
\arr{0.3 0.7}{0.7 0.3}
\arr{1.7 0.7}{1.3 0.3}
\arr{2.3 0.7}{2.7 0.3}

\arr{4.3 0.7}{4.7 0.3}
\arr{5.7 0.7}{5.3 0.3}
\arr{6.3 0.7}{6.7 0.3}
\put{$\ssize \gamma_1$} at 0.7 0.7
\put{$\ssize \gamma_2$} at 1.3 0.7
\put{$\ssize \gamma_{r-1}$} at 5.2 0.7
\put{$\ssize \gamma_r$} at 6.73 0.7
\endpicture} at 9 0
\endpicture}
$$
In all cases, the vertex $a(r)$ is a sink.
       \medskip
{\it The module $P_r(p)$ exists} (since for $r$ even, $p$ is a sink, whereas
for $r$ odd, $p$ is a source). Also, {\it the module $P_0(q)$ exists}
(since $q = a(r)$ is a sink) {\it and it is a submodule of $P_r(p)$.}

       \medskip
{\bf I. Case $b = r$.} We start with the path $\pi = (p|\gamma_1,\dots,\gamma_r|q)$
and obtain the module $P_r(p)$ and its
submodule $P_0(q)$. We define $M_0 = P_r(p)/P_0(q).$
This is obviously an incomplete sink module and we let $M_i = \sigma^i M_0$. 

First, let us consider $i = 1.$ We apply $\sigma$ to the
exact sequence
$$
  0 \to P_0(q) \to P_r(p) @>\epsilon>> M_0 \to 0
$$
  and obtain an exact sequence
$$
    0 \to P_{r-1}(p) @>\sigma\epsilon>> \sigma M_0 @>>> I_0(q) \to 0, \tag{$*$}
$$
since $\sigma P_{r}(p) = P_{r-1}(p)$. The support of $P_{r-1}(p)$
are the vertices $a$ of $T(n)$ with distance $d(p,a) \le r-1$.
In particular, the vertices $a(i)$ with $0\le i < r$ belong to this support.
On the other hand, the support of $I_0(q)$ is just the vertex $q$.
The arrow $\gamma_r$ is the only connection between the support of
$P_{r-1}(p)$ and the support of $I_0(q)$. In particular, we have
$\dim \Ext^1(I_0(q), P_{r-1}(p)) = 1$. This shows that $\sigma M_0$ is the
unique indecomposable extension of $P_{r-1}(p)$ by $I_0(q)$; it is given
by replacing in the direct sum $P_{r-1}(p)\oplus I_0(q)$ the zero map
at the arrow $\gamma_r$ by the identity map
$$
      (P_{r-1}(p))_{a(r-1)} = k @<1<< k = (I_0(q))_{a(r)}.
$$

Next, we apply $\sigma^{i-1}$ to $(*)$  with $1\le i-1 \le r-1$ and we
obtain the exact sequence
$$
  0 \to P_{r-i}(p) @>\sigma\epsilon>> \sigma^{i} M_0 @>>> I_{i-1}(q) \to 0.
$$
The arrow $\gamma_{r-i+1}$ is the only connection between the support of
$P_{r-i}(p)$ and the support of $I_{i-1}(q)$. In particular, we have
$\dim \Ext^1(I_{i-1}(q), P_{r-i}(p)) = 1$.
$$
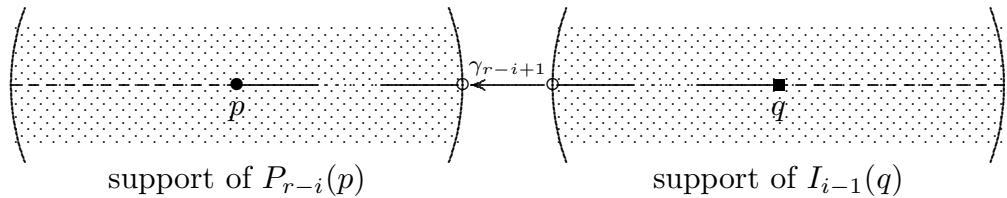
\hbox{\beginpicture
  \setcoordinatesystem units <1cm,.7cm>
  \multiput{} at 0 1  0 -1 /
\put{$\bullet$} at 0 0
\multiput{$\circ$} at  2.5 0  3.5 0 /
\linethickness1.25mm
\putrule from 5.94 0 to 6.06 0

\arr{3.4 0}{2.6 0}

\put{$\ssize\gamma_{r-i+1}\strut$} at 3 0.35

\plot 0 0  .9 0 /
\plot 1.6 0  2.45 0 /
\plot 3.55 0  4.4 0 /
\plot 5.1 0  6 0 /
\setdots <1mm>
\plot 1 0  1.5 0 /
\plot 4.5 0  5 0 /
\setsolid
\circulararc 20 degrees from 2.5 0 center at 0 0
\circulararc -20 degrees from 2.5 0 center at 0 0
\circulararc 20 degrees from -2.5 0 center at 0 0
\circulararc -20 degrees from -2.5 0 center at 0 0

\circulararc 20 degrees from 3.5 0 center at 6 0
\circulararc -20 degrees from 3.5 0 center at 6 0
\circulararc 20 degrees from 8.5 0 center at 6 0
\circulararc -20 degrees from 8.5 0 center at 6 0

\setdashes <1mm>
\plot 0 0  -2.5 0 /
\plot 6 0  8.5 0 /

\put{support of $P_{r-i}(p)$} at 0 -1.5
\put{support of $I_{i-1}(q)$} at 6 -1.5
\setshadegrid span <.7mm>
\vshade -2.5 -1 1  2.5 -1 1  /
\vshade 3.5 -1 1  8.5 -1 1  /
\put{$q\strut$} at 6 -.35
\put{$p\strut$} at 0 -.35
\endpicture}
$$
and $\sigma^i M_0$ is the
unique indecomposable extension of $P_{r-i}(p)$ by $I_{i-1}(q)$; it is given
by replacing in the direct sum $P_{r-i}(p)\oplus I_{i-1}(q)$ the zero map
at the arrow $\gamma_{r-i+1}$ by the identity map.

Finally, we apply $\sigma$ to the exact sequence
$$
  0 \to P_{0}(p) @>>> M_r @>f>> I_{r-1}(q) \to 0
$$
and obtain the exact sequence
$$
    0 @>>> M_{r+1} @>\sigma f>> I_r(q) \to I_0(p) \to 0,
$$
in some sense dual to the exact sequence we started with.
Of course, we see in this way that $M_{r+1}$ is an incomplete source module.
	\medskip 
Let us add how the modules $M_{-i}$ and $M_{r+1+i}$ with $i\ge 0$ look.

First, we consider the modules $M_{-i}$ with $i\ge 0$. They are given by a
projective presentation of the form
$$
 0 \to P_i(q) \to P_{r+i}(p) \to M_{-i} \to 0;
$$
since $\dim\Hom(P_i(q),P_{r+i}(p)) = \dim\Hom(P_0(q),P_r(p)) = P_r(p)_q = 1$,
the module $M_{-i}$ is uniquely defined in this way.
Clearly, the support of $M_{-i}$ is contained in the support of $P_{r+i}(p)$, and
this is the ball $B_{r+i}(p)$ with center $p$ and radius $r+i$. If
$$
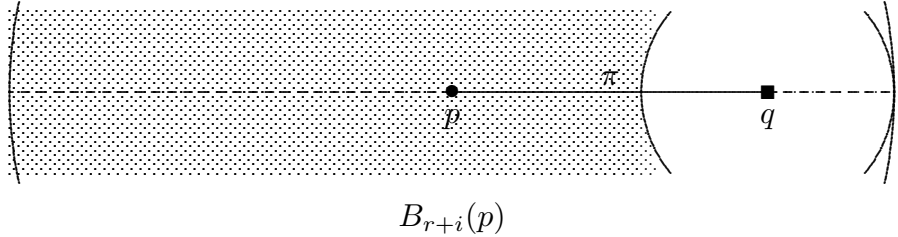
\hbox{\beginpicture
  \setcoordinatesystem units <.7cm,.7cm>
  \multiput{} at 2 1.5  -12 -1.5 /
\put{$\bullet$} at -5 0
\put{$\sssize\blacksquare$} at 0 0

\put{$\pi\strut$} at -2.5 0.25
\put{$p\strut$} at -5 -.35
\put{$q\strut$} at 0 -.35
\plot 0 0  -5 0 /
\setdots <1mm>
\plot 2 0  0 0 /
\plot -5 0  -12 0 /
\setsolid
\circulararc 12 degrees from 2 0 center at -5 0
\circulararc -12 degrees from 2 0 center at -5 0
\circulararc 12 degrees from -12 0 center at -5 0
\circulararc -12 degrees from -12 0 center at -5 0

\circulararc 40 degrees from 2 0 center at 0 0
\circulararc -40 degrees from 2 0 center at 0 0
\circulararc 40 degrees from -2 0 center at 0 0
\circulararc -40 degrees from -2 0 center at 0 0

\setdashes <.9mm>
\plot 0 0  2 0 /
\plot -5 0  -12 0 /

\put{$B_{r+i}(p)$} at -5 -2
\put{} at -5 -2
\setshadegrid span <.5mm>

\vshade -12 -1.3 1.3   <z,z,,> -2.06 -1.3 1.3  <z,z,,> -2.05 -1.3 0 <z,z,,> -1.75  -1.3 -1.3 /
\vshade -12 -1.3 1.3   <z,z,,> -2.06 -1.3 1.3  <z,z,,> -2.05 0 1.3 <z,z,,> -1.75  1.3 1.3 /
\endpicture}
$$
Any path in $B_{r+i}$ starting at $p$, ending in a leaf and not using $\gamma_1$ 
belongs to the support of $M_{-i}$. Since $n\ge 3$, we obtain in this way
diameter paths for $M_{-i}$ of length $2(r+i)$. This shows that $M_{-i}$ is a sink
module. 

Dually, the module $M_{r+1+i}$ with $i\ge 0$ has an injective copresentation of the form
$$
 0 \to M_{r+1+i} \to I_{r+i}(q) \to I_{i}(p) \to 0,
$$
thus its support is contained in the support of $I_{r+i}(q)$ 
and this is again is a ball with radius $r+i$, but now with center $q$.
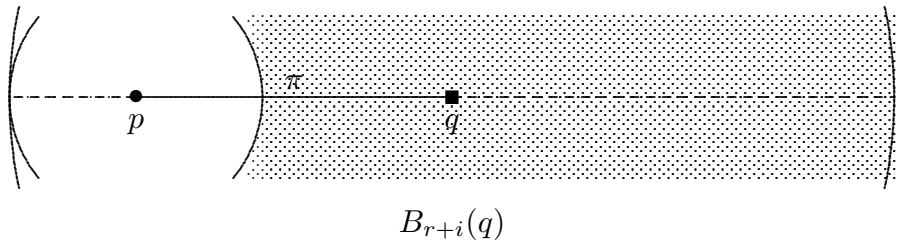
$$
\hbox{\beginpicture
  \setcoordinatesystem units <.7cm,.7cm>
  \multiput{} at -2 1.5  12 -1.5 /
\put{$\sssize\blacksquare$} at 5 0
\put{$\bullet$} at 0 0

\put{$\pi\strut$} at 2.5 0.25
\put{$q\strut$} at 5 -.35
\put{$p\strut$} at 0 -.35
\plot 0 0  5 0 /
\setdots <1mm>
\plot -2 0  0 0 /
\plot 5 0  12 0 /
\setsolid
\circulararc 12 degrees from -2 0 center at 5 0
\circulararc -12 degrees from -2 0 center at 5 0
\circulararc 12 degrees from 12 0 center at 5 0
\circulararc -12 degrees from 12 0 center at 5 0

\circulararc 40 degrees from -2 0 center at 0 0
\circulararc -40 degrees from -2 0 center at 0 0
\circulararc 40 degrees from 2 0 center at 0 0
\circulararc -40 degrees from 2 0 center at 0 0

\setdashes <.9mm>
\plot 0 0  -2 0 /
\plot 5 0  12 0 /

\put{$B_{r+i}(q)$} at 5 -2
\put{} at 0 -2
\setshadegrid span <.5mm>
\vshade 1.75  -1.3 -1.3 <z,z,,> 2.05 -1.3 0 <z,z,,> 2.06 -1.3 1.3 <z,z,,> 12 -1.3 1.3  /
\vshade 1.75  1.3 1.3 <z,z,,> 2.05 0 1.3  <z,z,,> 2.06 -1.3 1.3 <z,z,,> 12 -1.3 1.3  /

\endpicture}
$$
Clearly, $M_{r+1+i}$ is a source module.
	\bigskip
For $b = r = 4,$ and $-1\le i \le 6$, the modules $M_i$ have the following shapes:
$$
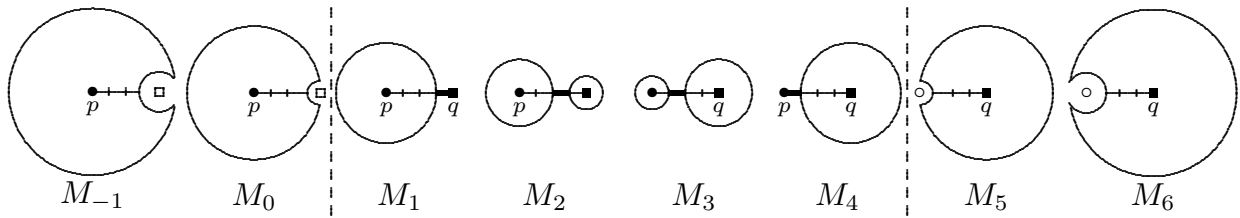
\hbox{\beginpicture
  \setcoordinatesystem units <.185cm,.185cm>
\multiput{} at 0 9  0 -3 /
\put{\beginpicture
\put{$M_{-1}$} at 0 -6.25
\plot 0 0  2.8 0 /
\put{$\ssize \bullet$} at 0 0
\plot 1 .2  1 -.2 /
\plot 2 .2  2 -.2 /
\circulararc 140 degrees from 2.8 0  center at 4 0
\circulararc -140 degrees from 2.8 0  center at 4 0
\circulararc 170 degrees from -5 0  center at 0 0
\circulararc -170 degrees from -5 0  center at 0 0
\plot 3.75 0.25  4.25 0.25  4.25 -.25  3.75 -.25  3.75 .25 /
\put{$\ssize p$} at 0 -1
\endpicture} at -19.5 0
\put{\beginpicture
\put{} at 4 0 
\put{$M_0$} at 0 -6.2
\plot 0 0  3.2 0 /
\put{$\ssize \bullet$} at 0 0
\plot 1 .2  1 -.2 /
\plot 2 .2  2 -.2 /
\circulararc 92 degrees from 3.2 0  center at 3.9 0
\circulararc -92 degrees from 3.2 0  center at 3.9 0
\circulararc 170 degrees from -4 0  center at 0 0
\circulararc -170 degrees from -4 0  center at 0 0
\plot 3.75 0.25  4.25 0.25  4.25 -.25  3.75 -.25  3.75 .25 /
\put{$\ssize p$} at 0 -1

\endpicture} at -8.5 0
\put{\beginpicture
\put{$M_1$} at .8 -6.2
\plot 0 0  4 0 /
\put{$\ssize \bullet$} at 0 0
\plot 1 .2  1 -.2 /
\plot 2 .2  2 -.2 /
\plot 3 .2  3 -.2 /
\circulararc 360 degrees from 3 0  center at 0 0
\linethickness1mm
\putrule from 3.75 0 to 4.25 0
\linethickness.5mm
\putrule from 3 0 to 4 0
\put{$\ssize q$} at 4 -1
\put{$\ssize p$} at 0 -1
\endpicture} at 0 0
\put{\beginpicture
\put{$M_2$} at 1.5 -6.2
\plot 0 0  4 0 /
\put{$\ssize \bullet$} at 0 0
\plot 1 .2  1 -.2 /
\plot 2 .2  2 -.2 /
\plot 3 .2  3 -.2 /
\circulararc 360 degrees from 2 0  center at 0 0
\circulararc 360 degrees from 3 0  center at 4 0
\linethickness1mm
\putrule from 3.75 0 to 4.25 0
\linethickness.5mm
\putrule from 2 0 to 3 0
\put{$\ssize p$} at 0 -1
\endpicture} at 8 0
\put{\beginpicture
\put{$M_3$} at 2.5 -6.2
\plot 0 0  4 0 /
\put{$\ssize \bullet$} at 0 0
\plot 1 .2  1 -.2 /
\plot 2 .2  2 -.2 /
\plot 3 .2  3 -.2 /
\circulararc 360 degrees from 1 0  center at 0 0
\circulararc 360 degrees from 2 0  center at 4 0
\linethickness1mm
\putrule from 3.75 0 to 4.25 0
\linethickness.5mm
\putrule from 1 0 to 2 0
\put{$\ssize q$} at 4 -1
\endpicture} at 16 0
\put{\beginpicture
\put{$M_4$} at 3.2 -6.2
\plot 0 0  4 0 /
\put{$\ssize \bullet$} at 0 0
\plot 1 .2  1 -.2 /
\plot 2 .2  2 -.2 /
\plot 3 .2  3 -.2 /
\circulararc 360 degrees from 1 0  center at 4 0
\linethickness1mm
\putrule from 3.75 0 to 4.25 0
\linethickness.5mm
\putrule from 0 0 to 1 0
\put{$\ssize q$} at 4 -1
\put{$\ssize p$} at 0 -1
\endpicture} at 24 0
\put{\beginpicture
\put{$M_5$} at 4 -6.2
\plot 0.8 0  4 0 /
\put{$\ssize \circ$} at 0 0
\plot 2 .2  2 -.2 /
\plot 3 .2  3 -.2 /
\circulararc 92 degrees from 0.8 0  center at 0.1 0
\circulararc -92 degrees from 0.8 0  center at 0.1 0
\circulararc 170 degrees from 8 0  center at 4 0
\circulararc -170 degrees from 8 0  center at 4 0
\linethickness1mm
\putrule from 3.75 0 to 4.25 0
\put{$\ssize q$} at 4 -1

\endpicture} at 32.5 0
\put{\beginpicture
\put{$M_6$} at 4 -6.2
\plot 1.2 0  4 0 /
\put{$\ssize \circ$} at 0 0
\plot 2 .2  2 -.2 /
\plot 3 .2  3 -.2 /
\circulararc 150 degrees from 1.2 0  center at 0 0
\circulararc -150 degrees from 1.2 0  center at 0 0
\circulararc 170 degrees from 9 0  center at 4 0
\circulararc -170 degrees from 9 0  center at 4 0
\linethickness1mm
\putrule from 3.75 0 to 4.25 0
\put{$\ssize q$} at 4 -1
\endpicture} at 42.5 0
\setdashes <1mm>
\plot -5.2 -4  -5.2 8.5 /
\plot 29.3 -4  29.3 8.5 /
\endpicture}
$$
Note that for $i$ even, the vertex $p$ (drawn as a bullet or a small circle) 
is a sink for $T(M_i)$ and $q$ (drawn as a small square) is a source,
whereas for $i$ odd, $p$ is a source and $q$ a sink. 

    \bigskip

{\bf II. Case $b < r$, with $(n,b,r)\neq (3,0,1)$.}
      \medskip
Let $b < r$ and $(n,b,r)\neq (3,0,1)$.
Here is the recipe for $M_0$. Let $(p|\gamma_1,\dots,\gamma_b|q)$ be a path of length $b$
and 
$
  (q|\delta_{b+1},\dots,\delta_r|x),\   (q|\delta'_{b+1},\dots,\delta'_r|x')
$
paths of length $r\!-\!b$ with $\delta_{b+1} \neq \delta'_{b+1}$:
$$
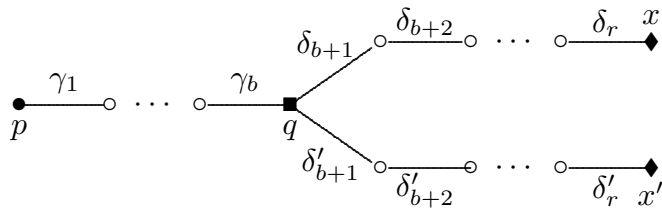
\hbox{\beginpicture
  \setcoordinatesystem units <1cm,.7cm>
\put{$\bullet$} at  0 0
\put{$\sssize \blacksquare$} at 3 0
\multiput{$\ssize \blacklozenge$} at 7 1  7 -1 /
\multiput{$\circ$} at 1 0  2 0  4 1  4 -1  5 1  5 -1  6 1  6 -1 /
\plot 0 0  0.93 0 /
\plot 2.1 0  3 0  3.9 .9 /
\plot 4.1 1  4.9 1 /
\plot 6.1 1  7 1 /
\plot 6.1 -1  7 -1 /
\plot 3 0  3.9 -.9 /
\plot 4.1 -1  5 -1 /
\put{$p$} at 0 -.4
\put{$q$} at 3 -.4
\put{$x$} at 7 1.4
\put{$x'$} at 7 -1.4
\put{$\gamma_1$} at .5 0.35
\put{$\gamma_b$} at 2.5 0.35
\multiput{$\cdots$} at 1.5 0  5.5 1  5.5 -1 /
\put{$\delta_{b+1}$} at 3.4 0.95
\put{$\delta'_{b+1}$} at 3.45 -.95
\put{$\delta_{b+2}$} at 4.5 1.3
\put{$\delta'_{b+2}$} at 4.5 -1.35
\put{$\delta_{r}$} at 6.5 1.3
\put{$\delta'_{r}$} at 6.5 -1.35
\endpicture}
$$
Let
$$
 M_0 = P_r(p)/(P_0(x)\oplus P_0(x')).
$$
As we will see, this is the incomplete sink module we are looking for.
(The module $M_0$ can be defined also in case $(n,b,r) = (3,0,1)$, but then $M_0$ is a flow
module.) 
      \medskip
We claim that $M_0$ is a module with center $p$ and diameter $2r$
and that $M_{b+1} = \sigma^{b+1}M_0$ has center $q$ and diameter $2r$.
This then implies that the modules $M_1,\dots,M_b$ have to be
flow modules and that $M_0$ is the incomplete sink module,
$M_{b+1}$ the incomplete source module in the $\sigma$-orbit
(see the Corollary at the end of section 3).
       \medskip
In order to see that $M_0$ is a (sink) module with center $p$ and diameter $2r$, we
have to exhibit a diameter path of the support
$B_r(p)$ of $P_r(p)$
which does not start or end in $x$  or $x'$. Then this is a diameter path for $M_0$,
has center $p$ and length $2r$.
In case $n\ge 4$, or $b\ge 1$, we take a diameter path of $B_r(p)$ which does not involve
$\gamma_1$. If
$b = 0$ and $r\ge 2$, let $z,z'$ be leaves of $B_r(p)$ with $d(x,z) = 2$ and
$d(x',z') = 2$. Then the path of $P_r(p)$ which connects $z$ with $z'$ goes through $p$,
thus is a diameter path of $M_0$.

Now let us consider the modules $M_i = \sigma^iM_0$ with $1\le i \le b+1.$ The
defining exact sequence
$$
 0 \to P_0(x)\oplus P_0(x') \to P_r(p) \to M_0 \to 0 
$$
yields an exact sequence
$$
 0 \to P_{r-1}(p) \to M_1 \to I_0(x)\oplus I_0(x') \to 0. \tag{$*$}
$$
If we apply $\sigma^b$ to this exact sequence $(*)$, we obtain a corresponding exact sequence
$$
 0 \to P_{r-b-1}(p) \to M_{b+1} \to I_b(x)\oplus I_b(x') \to 0
$$
(here we use that $b < r$). This shows that the support of $M_{b+1}$ is the
union of the ball $B_{r-b-1}(p)$ and the balls $B_b(x)$ and $B_b(x').$

First of all, let us construct a path of length $2r$ in the support of $M_{b+1}$.
We choose a path in $B_b(x)$ from the center $x'$ to a boundary vertex, say $w$,
not using the arrow
$\delta_r$. Similarly, we choose a path in $B_b(x')$
from the center $x$ to a boundary vertex, say $w'$,
not using the arrow $\delta'_r$. Combining these two paths with the given path from $x$ via $q$
to $x'$ (using the arrows $\delta_i$ and $\delta'_i$)
we obtain a path of length $2b + 2(r-b) = 2r$ with center $q$.

On the other hand, we claim that $T(M_{b+1}) \subseteq B_r(q).$
Namely, if $u$ belongs to $B_b(x)$,
then $d(u,q) \le d(u,x)+d(x,q) \le b+(r-b) = r.$ Similarly, for $u\in B_b(x')$, we have
$d(u,q)\le r.$ Finally, if $u\in B_{r-b-1},$ then $d(u,q) \le d(u,p)+d(p,q) \le
(r-b-1)+b = r-1.$  This completes the proof. Actually, the last calculations show that
any boundary vertex of $T(M_{b+1})$ is a boundary vertex of $B_b(x)$ or of $B_b(x')$
(thus, in particular, a source).
   \bigskip

It seems to be of interest to look at the module $M_{b+1}$ in more detail.
We distinguish two cases, the case $r> 2b$ and the case $r\le 2b$.

First, let us assume that $r > 2b$, thus $b < r-b.$ In this case, the extension
of the module $P_{r-b-1}(p)$ by the module $I_b(x)$ is furnished by the arrow
$\delta_{r-b}$ (and similarly, the extension
of the module $P_{r-b-1}(p)$ by the module $I_b(x')$ is furnished by the arrow
$\delta'_{r-b}$). In particular, the three balls $B_{r-b-1}(p), B_b(x)$ and $B_b(x')$
are pairwise disjoint.

On the other hand, if $r\le 2b$, thus $r-b \le b$, the extension of $P_{r-b-1}(p)$
by both modules $I_b(x)$ and $I_b(x')$ is furnished by the arrow
$\gamma_{r-b}$. In this case, the vertex $q$ belongs both to $B_b(x)$ and $B_b(x').$
	\medskip

In general, the modules $M_i$ with $1\le i \le r$ are quite easy to visualize. We apply
$\sigma^{i-1}$ to the exact sequence $(*)$ and obtain the exact sequence:
$$
 0 \to P_{r-i}(p) \to M_{i} \to I_{i-1}(x)\oplus I_{i-1}(x') \to 0.
$$

First, assume that $i \le r-b$. Then the extension of $P_{r-i}(p)$ by $I_{i-1}(x)$  or
by $I_{i-1}(x')$ is given by $\delta_{r-i+1}$ or $\delta'_{r-i+1}$, respectively:
$$
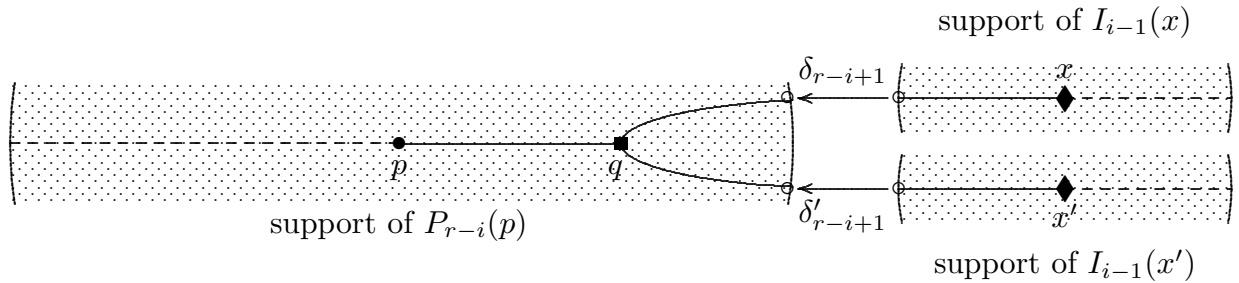
\hbox{\beginpicture
  \setcoordinatesystem units <1.23cm,.5cm>
  \multiput{} at -3.5 1.5  7.5 -1.5 /
\put{$\bullet$} at 0 0 
\plot  0 0  2 0 /
\linethickness1.25mm
\putrule from 1.94 0 to 2.06 0
\circulararc 9 degrees from -3.5 0 center at 0 0
\circulararc -9 degrees from -3.5 0 center at 0 0
\circulararc 9 degrees from 3.55 0 center at 0 0
\circulararc -9 degrees from 3.55 0 center at 0 0
\setdashes <1mm>
\plot -3.5 0  0 0 /
\setsolid
\ellipticalarc axes ratio 5:1 77 degrees from 2 0  center at 4 0
\ellipticalarc axes ratio 5:1 -77 degrees from 2 0  center at 4 0

\multiput{$\circ$} at  3.5 1  4.5 1 /
\multiput{$\circ$} at  3.5 -1  4.5 -1 /
\arr{4.4 1}{3.6 1}
\arr{4.4 -1}{3.6 -1}
\put{$\delta_{r-i+1}$} at 4 1.6
\put{$\delta'_{r-i+1}$} at 4 -1.6
\multiput{$\blacklozenge$} at 6 1  6 -1 /
\plot 4.5 1  6 1 /
\plot 4.5 -1  6 -1 /
\setdashes <1mm>
\plot 6 1  7.5 1 /
\plot 6 -1  7.5 -1 /
\setsolid
\circulararc 12 degrees from 4.5 1 center at 6 1
\circulararc -12 degrees from 4.5 1 center at 6 1

\circulararc 12 degrees from 4.5 -1 center at 6 -1
\circulararc -12 degrees from 4.5 -1 center at 6 -1

\circulararc 12 degrees from 7.5 1 center at 6 1
\circulararc -12 degrees from 7.5 1 center at 6 1

\circulararc 12 degrees from 7.5 -1 center at 6 -1
\circulararc -12 degrees from 7.5 -1 center at 6 -1
\put{$p$} at 0 -0.6
\put{$q$} at 1.95 -0.6
\put{$x$} at 6 1.6
\put{$x'$} at 6 -1.6
\setshadegrid span <.7mm>
\vshade -3.5 -1.4 1.4  3.5 -1.4 1.4  /

\vshade 4.5 .2 1.8  7.5 .2 1.8  /
\vshade 4.5 -1.8 -.2  7.5 -1.8 -.2  /
\put{support of $P_{r-i}(p)$} at 0 -1.8
\put{support of $I_{i-1}(x)$} at 6  2.7
\put{support of $I_{i-1}(x')$} at 6  -2.7
\endpicture}
$$

Second, assume that $i > r-b$. Then the extensions of $P_{r-i}(p)$ by $I_{i-1}(x)$  as well as 
by $I_{i-1}(x')$ are given by $\gamma_{r-i+1}$:
$$
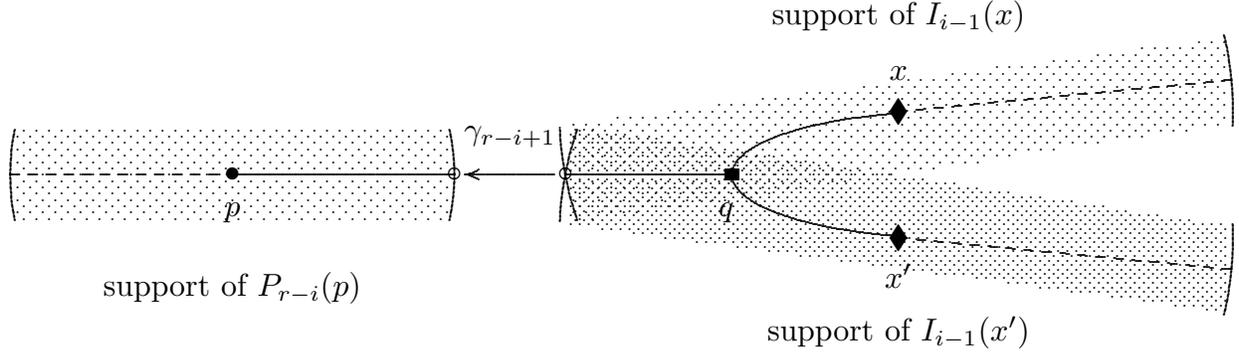
\hbox{\beginpicture
  \setcoordinatesystem units <1.23cm,.7cm>
  \multiput{} at -2 1.5  9 -1.5 /
\put{$\bullet$} at 0 0 
\plot  0 0  2 0 /
\multiput{$\circ$} at  2 0  3 0  /
\arr{2.9 0}{2.1 0}
\put{$\gamma_{r-i+1}$} at 2.5 0.6
\plot 3 0  4.5 0 /
\linethickness1.25mm
\putrule from 4.44 0 to 4.56 0
\circulararc 12 degrees from -2 0 center at 0 0
\circulararc -12 degrees from -2 0 center at 0 0
\circulararc 12 degrees from 2 0 center at 0 0
\circulararc -12 degrees from 2 0 center at 0 0
\setdashes <1mm>
\plot -2 0  0 0 /
\setsolid

\ellipticalarc axes ratio 3.5:1 77 degrees from 4.5 0  center at 6.5 0
\ellipticalarc axes ratio 3.5:1 -77 degrees from 4.5 0  center at 6.5 0

\multiput{$\blacklozenge$} at 6 1  6 -1 /
\setdashes <1mm>
\plot 6 1  9 1.5 /
\plot 6 -1  9 -1.5 /
\setsolid
\circulararc 8 degrees from 3 0 center at 6 1
\circulararc -8 degrees from 3 0 center at 6 1

\circulararc 8 degrees from 3 0 center at 6 -1
\circulararc -8 degrees from 3 0 center at 6 -1

\circulararc 8 degrees from 9 1.5 center at 6 1
\circulararc -8 degrees from 9 1.5 center at 6 1

\circulararc 8 degrees from 9 -1.5 center at 6 -1
\circulararc -8 degrees from 9 -1.5 center at 6 -1
\put{$p$} at 0 -0.6
\put{$q$} at 4.45 -0.6
\put{$x$} at 6 1.6
\put{$x'$} at 6 -1.6
\setshadegrid span <.7mm>
\vshade -2 -.8 .8  2 -.8 .8  /

\vshade 3 -.8 .8  9 .8 2.2  /
\setshadegrid span <.5mm>
\vshade 3 -.8 .8  9 -2.2 -.8  /
\put{support of $P_{r-i}(p)$} at 0 -1.8
\put{support of $I_{i-1}(x)$} at 6  2.5
\put{support of $I_{i-1}(x')$} at 6  -2.5
\endpicture}
$$
	\medskip
We also want to describe the module $M_{r+1}$. The exact sequence
$$
 0 \to P_{0}(p) \to M_{r} \to I_{r-1}(x)\oplus I_{r-1}(x') \to 0
$$
yields the exact sequence
$$
 0 \to M_{r+1} \to I_{r}(x)\oplus I_{r}(x') \to I_0(p) \to 0.
$$
This shows that $M_{r+1}$ is a maximal submodule of $I_{r}(x)\oplus I_{r}(x')$
with factor the simple module $I_0(p).$ 
       	     \bigskip\bigskip
{\bf III. Case $b = 0, r = 1$.} This concerns the remaining case $(n,b,r) = (3,0,1)$, but
works for all $n\ge 3$.
           \smallskip
We start with a source $p$ and take an indecomposable module $M$ with
$\dim M_p = 2$ and $\dim M_x = 1$ for all the neighbors
$x$ of $p$, whereas $M_a = 0$ for all other vertices $a$ of $T(n)$ (for $n= 3$, this
concerns the representations of a quiver of type $\Bbb D_4$, thus
there is a unique such module; for $n\ge 4$, there are many such modules).
Since $T(M)$ is the full subquiver of $T(n)$ given by $p$ and its neighbors,
the module $M$ is a sink module with radius $1$. Since $\dim M_x = 1 < 2 = \dim M_p$, for
the neighbors $x$ of $p$, we see that $M$ is incomplete.
The vertices of the support of $\sigma M$ are again $p$ (now a source)
and its neighbors $x$ (now sinks), and $\dim (\sigma M)_p = 2, \dim (\sigma M)_x = 1$.
We see that $M_1 = \sigma M$ is an incomplete source module with center $q = p$ and radius $1$.
The center path for the shift orbit of $M$ is the path
from $p$ to $q = p$ of length $b = 0$. \hfill$\square$
	\medskip
This completes the proof of Theorem 2. 
\hfill$\square$

	\bigskip
{\bf Further examples.}
The examples constructed until now were individual modules.
The same techniques allow to exhibit also families of modules. Let us construct 
a $1$-parameter family of incomplete sink modules $M_0$
with a fixed dimension vector and equal radius and center path.

We start with a subquiver of $T(n)$ of the following kind:
$$
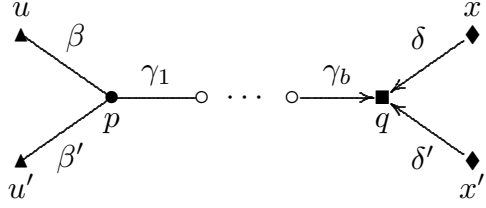
\hbox{\beginpicture
  \setcoordinatesystem units <1cm,.7cm>
  \put{$\bullet$} at  0 0
  \put{$\sssize \blacksquare$} at 3 0
  \multiput{$\ssize \blacklozenge$} at 4 1  4 -1 /
  \multiput{$\circ$} at 1 0  2 0   /
  \plot 0 0  0.93 0 /
  \arr{2.1 0}{2.9 0}  
  \arr{3.9 .9}{3.1 0.1}
  \arr{3.9 -.9}{3.1 -0.1}
  \put{$p$} at 0 -.4
  \put{$q$} at 3 -.4
  \put{$x$} at 4 1.4
  \put{$x'$} at 4 -1.4
  \put{$u$} at -1 1.4
  \put{$u'$} at -1 -1.4
  \put{$\gamma_1$} at .5 0.35
  \put{$\gamma_b$} at 2.5 0.35
  \multiput{$\cdots$} at 1.5 0  /
  \put{$\delta$} at 3.4 0.95
    \put{$\delta'$} at 3.45 -.95
\multiput{$\ssize \blacktriangle$} at -1 1  -1 -1 /
\plot -1 1  0 0  -1 -1 /
  \put{$\beta$} at -.4 0.95
    \put{$\beta'$} at -.45 -.95
\endpicture}
$$
(with $x,x'$ being sources). Let $\Cal A$ be the set of the following four
representations of $(T(n),\Omega)$
$$
  P_{b+1}(u),\ P_{b+1}(u'),\ I_0(x),\ I_0(x').
  $$ 
These are pairwise orthogonal bricks: one only has to verify that the first
two modules are orthogonal bricks, however $P_{b+1}(u) = \sigma^{-b-1}P_0(u)$,
$P_{b+1}(u') = \sigma^{-b-1}P_0(u')$, and $P_0(u),\ P_0(u')$ are of course
orthogonal bricks. Let us denote by $\Cal E(\Cal A)$ the extension closure of
$\Cal E$. As one knows from [R1], $\Cal E(\Cal A)$ is an exact abelian subcategory
with simple objects the modules $\Cal A$, and that this subcategory is equivalent
to the representations of the $\Ext$-quiver $\Delta(\Cal A)$ of $\Cal A$: the vertices
of $\Delta(\Cal A)$ are of the form $[A]$, where
$A$ is an object in $\Cal A$, and there are $t$ arrows $[A] \to [A']$ provided
$\dim \Ext^1(A,A') = t.$ In our case, $\Delta(\Cal A)$ is the following bipartite quiver
$$
\hbox{\beginpicture
  \setcoordinatesystem units <1.5cm,.9cm>
\multiput{$\circ$} at  0 0  1 0  0 1  1 1 /
\arr{0.9 0.9}{0.1 0.1}
\arr{0.9 0}{0.1 0}
\arr{0.9 0.1}{0.1 0.9}
\arr{0.9 1}{0.1 1}
\put{$P_{b+1}(u)$} at -0.5 1
\put{$P_{b+1}(u')$} at -0.5 0
\put{$I_0(x)$} at 1.4 1
\put{$I_0(x')$} at 1.4 0

\endpicture}
$$
thus an affine quiver of type $\widetilde {\Bbb A}_{22}$.

Let $\Cal E'(\Cal A)$ be the class of indecomposable objects in $\Cal E(\Cal A)$
which correspond to sincere representations of the quiver $\Delta(\Cal A)$,
thus they are the indecomposable representations $M$ of $(T(n),\Omega)$
with an exact sequence of the form
$$
 0 \to    P_{b+1}(u)^s\oplus P_{b+1}(u')^{s'} \to M \to I_0(x)^t\oplus I_0(x')^{t'}\to 0 \tag{$*$}
 $$
 with positive integers $s,s',t,t'$.

We claim that $M$ is a sink module with center $p$ (and radius $b+2$) and that
$\sigma^{b+1}M$ is a source module with center $q$ (and the same radius $b+2$).
Since $d(p,q) = b$,
it follows that $M$ is an incomplete sink module
with center path $(p|\gamma_1,\dots,\gamma_b|q)$ (see the Corollary at the end of section 3).

Thus, let us analyze the diameter paths of $M$. It is obvious that $T(M) \subseteq B_{b+2}(p)$.
Since $s\ge 1$, there is a path from the boundary of $B_{b+1}(u)$ to its center $u$
which does not use $\beta$. Similarly, since $s'\ge 1$, there is a path from the
boundary of $B_{b+1}(u')$ to its center $u'$ which does not use $\beta'$. Combining these two
paths with the arrows $\beta$ and $\beta'$, we obtain a path of length $2b+4,$ thus
a diameter path and its center is $p$.
This shows that $M$ is a sink module with center $p$ and radius $b+2$.

If we apply $\sigma^{b+1}$ to $(*)$, we obtain the exact sequence
$$
 0 \to P_{0}(u)^s\oplus P_{0}(u')^{s'} \to \sigma^{b+1}M
  \to I_{b+1}(x)^t\oplus I_{b+1}(x')^{t'}\to 0.
$$
Thus, by duality we see that
$\sigma^{b+1}M$ is a source module with center $q$ and radius $b+2$.

Let us add that the flow modules in the $\sigma$-orbit of $M$ are the
modules $\sigma^iM$ with $1\le i \le b$, they are middle terms of exact
sequences of the following form:
$$
 0 \to P_{b+1-i}(u)^s\oplus P_{b+1-i}(u')^{s'} \to \sigma^{i}M
  \to I_{i}(x)^t\oplus I_{i}(x')^{t'}\to 0
  $$
and the extension is furnished by the arrow $\gamma_{b+1-i}.$

     \bigskip\bigskip
{\bf 6. The $\tau$-orbits.}
     \medskip
As we have mentioned in (2.6),
the Auslander-Reiten translation $\tau$ in $\mo(T(n),\Omega)$ is nothing else
than $\sigma^2$. Thus, Theorem 1 provides information on the $\tau$-orbits of the regular
indecomposable modules and on the shape of the components of the Auslander-Reiten quiver.

One should be aware that the very lucid behavior of the shift orbits looks more
complicated when we deal with the $\tau$-orbits. We use again the labeling $M_i$ of the
modules in a shift orbit, with $M_0$ the incomplete sink module and $M_i = \sigma^i M_0$.
As we have mentioned in the introduction, the modules of the form $M_i$ with $i$ 
even will be called {\it even} modules, those with
$i$ odd will be called {\it odd} modules.
 
We obtain in this way two $\tau$-orbits (but remember: in different categories --- one
$\tau$-orbit consists of representations of $(T(n),\Omega)$, the other of representations of
$(T(n),\sigma\Omega)$). Also, we have to stress that the labeling ``even'' and ``odd''
refers to our interest in sink modules. 
In a similar way, one may focus the attention to the unique incomplete source module in a
given $\sigma$-orbit. Alternatively, we may concentrate on the 
invariant $b$ which may be even or odd. It turns out that 
there are four different kinds of $\tau$-orbits of regular indecomposable module.
Here are these $\tau$-orbits (below any module, we mention its radius).
	\medskip
{\bf The even modules, with $b$ even:}
$$
\hbox{\beginpicture
  \setcoordinatesystem units <1cm,.6cm>
\put{$\cdots$} at -1 0
\put{$M_{-2}$} at 0 0
\put{$M_{0}$} at 2 0
\put{$M_{2}$} at 4 0
\put{$\cdots$} at 5 0
\put{$M_{b}$} at 7 0
\put{$M_{b+2}$} at 9 0
\put{$M_{b+4}$} at 11 0
\put{$\cdots$} at 12 0
\setdashes <1mm>
\plot 2.5 1.1 2.5 -1.3 /
\plot 7.5 1.1 7.5 -1.3 /
\put{sink modules} [r] at 1.5 1
\put{flow modules} at 5 1
\put{source modules} [l] at 8.5 1

\put{$\cdots$} at -1 -1
\put{$\ssize r+2$} at 0 -1
\put{$\ssize r$} at 2 -1
\put{$\ssize r\minus 1$} at 4 -1
\put{$\cdots$} at 5 -1
\put{$\ssize r\minus 1$} at 7 -1
\put{$\ssize r+1$} at 9 -1
\put{$\ssize r+3$} at 11 -1
\put{$\cdots$} at 12 -1

\put{} at 0 -2.5
\endpicture}
$$

{\bf The even modules with $b$ odd:}
$$
\hbox{\beginpicture
  \setcoordinatesystem units <1cm,.6cm>
\put{$\cdots$} at -1 0
\put{$M_{-2}$} at 0 0
\put{$M_{0}$} at 2 0
\put{$M_{2}$} at 4 0
\put{$\cdots$} at 5 0
\put{$M_{b-1}$} at 6 0
\put{$M_{b+1}$} at 8 0
\put{$M_{b+3}$} at 10 0
\put{$\cdots$} at 11 0
\setdashes <1mm>
\plot 2.5 1.1 2.5 -1.3 /
\plot 7.5 1.1 7.5 -1.3 /
\put{sink modules} [r] at 1.5 1
\put{flow modules} at 5 1
\put{source modules} [l] at 8.5 1

\put{$\cdots$} at -1 -1
\put{$\ssize r+2$} at 0 -1
\put{$\ssize r$} at 2 -1
\put{$\ssize r\minus 1$} at 4 -1
\put{$\cdots$} at 5 -1
\put{$\ssize r\minus 1$} at 6 -1
\put{$\ssize r$} at 8 -1
\put{$\ssize r+2$} at 10 -1
\put{$\cdots$} at 11 -1

\put{} at 0 -2.5
\endpicture}
$$

{\bf The odd modules with $b$ even:}
$$
\hbox{\beginpicture
  \setcoordinatesystem units <1cm,.6cm>
\put{$\cdots$} at -1 0
\put{$M_{-3}$} at 0 0
\put{$M_{-1}$} at 2 0
\put{$M_{1}$} at 4 0
\put{$\cdots$} at 5 0
\put{$M_{b-1}$} at 6 0
\put{$M_{b+1}$} at 8 0
\put{$M_{b+3}$} at 10 0
\put{$\cdots$} at 11 0
\setdashes <1mm>
\plot 3.5 1.1 3.5 -1.3 /
\plot 7.5 1.1 7.5 -1.3 /
\put{sink modules} [r] at 1.5 1
\put{flow modules} at 5 1
\put{source modules} [l] at 8.5 1

\put{$\cdots$} at -1 -1
\put{$\ssize r+3$} at 0 -1
\put{$\ssize r+1$} at 2 -1
\put{$\ssize r\minus 1$} at 4 -1
\put{$\cdots$} at 5 -1
\put{$\ssize r\minus 1$} at 6 -1
\put{$\ssize r$} at 8 -1
\put{$\ssize r+2$} at 10 -1
\put{$\cdots$} at 11 -1

\put{} at 0 -2.5
\endpicture}
$$

{\bf The odd modules with $b$ odd:}
$$
\hbox{\beginpicture
  \setcoordinatesystem units <1cm,.6cm>
\put{$\cdots$} at -1 0
\put{$M_{-3}$} at 0 0
\put{$M_{-1}$} at 2 0
\put{$M_{1}$} at 4 0
\put{$\cdots$} at 5 0
\put{$M_{b}$} at 6 0
\put{$M_{b+2}$} at 8 0
\put{$M_{b+4}$} at 10 0
\put{$\cdots$} at 11 0
\setdashes <1mm>
\plot 3.5 1.1  3.5 -1.3 /
\plot 6.5 1.1  6.5 -1.3 /
\put{sink modules} [r] at 1.5 1
\put{flow modules} at 5 1
\put{source modules} [l] at 8.5 1

\put{$\cdots$} at -1 -1
\put{$\ssize r+3$} at 0 -1
\put{$\ssize r+1$} at 2 -1
\put{$\ssize r\minus 1$} at 4 -1
\put{$\cdots$} at 5 -1
\put{$\ssize r\minus 1$} at 6 -1
\put{$\ssize r+1$} at 8 -1
\put{$\ssize r+3$} at 10 -1
\put{$\cdots$} at 11 -1

\put{} at 0 -2.5
\endpicture}
$$

As we have seen in section 3, looking at a shift orbit, the change of the radius 
of the corresponding modules is given by a simple and uniform rule, 
in contrast to the four different 
rules which occur for $\tau$-orbits.

	\bigskip
{\bf Proposition 2.} {\it Let $\Cal X$ be a $\tau$-orbit. 
We assume that the 
sink modules in $\Cal X$ have center $p$, the source modules in $\Cal X$ have center $q$
and that $b = d(p,q)$.}

(a) {\it If $b = 2s$ is even, then $\Cal X$ contains precisely $s$ flow modules.
}

(b) {\it If $b = 2s+1$ is odd, then the number of flow modules in $\Cal X$ 
is  $s$ or $s+1$. Let $X$ be a sink module in $\Cal X$ and that $\tau X$ is not a sink module.
Then $\iota(X) = 0$ in case the number of flow modules in $\Cal X$ is $s$, otherwise
$\iota(X) = -1$.}
\hfill$\square$
	\medskip
Of course, always we know: If $X$ is a sink module and $\tau X$ is not a sink module,
then $\iota(X)$ is equal to $-1$ or $0$. 
	\medskip
{\bf Proposition 3.} {\it Let $\Cal X$ be a $\tau$-orbit. 
We assume that the 
sink modules in $\Cal X$ have center $p$, the source modules in $\Cal X$ have center $q$
and that $b = d(p,q).$ We assume that $b = 2s$ is even. Let $X$ be a sink module in $\Cal X$
such that $\tau X$ is not a sink module, thus $\iota(X)$ is equal to $0$ or
$-1$. If  $\iota(X) = 0$, then 
$r(X) < r(\tau^{s+1} X)$. If  $\iota(X) = -1$, then 
$r(X) > r(\tau^{s+1} X)$.} 
	\medskip
Proof. Let $X_0$ be the incomplete sink module in the $\sigma$-orbit of $X$.
and $r_0$ the radius of $X_0$. By Theorem 1, the radius of $X_{-1}$ is $r_0+1$, 
the flow modules $X_1,\dots, X_{2s}$ have radius 
$r_0-1$, the module $X_{2s+1}$ has radius $r_0$ and the module
$X_{2s+2}$ has radius $r_0+1$.
For $\iota(X) = 0$, we have $\tau^{s+1}X = X_{2s+2}$, thus 
$r(X) = r_0 < r_0+1 = r(\tau^{s+1}X).$
Similarly, for $\iota(X) = -1$, we have $X = X_{-1}$ and $\tau^{s+1}X = X_{2s+1}$, thus
$r(X) = r_0+1 > r_0 = r(X_{2s+1}).$
\hfill$\square$
	\bigskip
{\bf 7. Auslander-Reiten components (Proof of Theorem 3 and Theorem 4).}
	\medskip
It is well-known that the category $\mo (T(n),\Omega)$ has Auslander-Reiten sequences
and that regular components are of the form $\Bbb Z\Bbb A_\infty$, 
thus the indecomposable regular modules are quasi-serial (see [R]);
we denote by $\ql(X)$ the quasi-length of $X$.
We say that a regular component $\Cal C$ of $\mo (T(n),\Omega)$
is {\it even} or {\it odd}, provided the quasi-simple
modules in $\Cal C$ are even or odd, respectively.

Given an indecomposable regular module $X$, we define
$$
 \overline\iota(X) = \frac1{\ql(X)} \sum\nolimits_F \iota(F)
$$
where we sum over all quasi-composition factors $F$ of $X$
and we call $\overline\iota(X)$ the {\it average index} of $X$. We will see that $\overline\iota(X) =
\iota(X).$ 
	\bigskip 
{\bf Lemma 3.} {\it Let $X$ be indecomposable and regular with $\overline\iota(X) \le 0$.
Then $X$ is a sink module and $B(X) = B(X')$, where $X'$ is the quasi-top of $X$.}
	\medskip 
Proof. We assume that the quasi-composition 
factors $F$ of $X$ are
$$
 M_i, M_{i+2},\dots, M_{i+2(l-1)},
$$
where $l$ is the quasi-length of $X$ (and where $M_0$ is an incomplete sink module
and $M_i = \sigma^i M_0$). We assume that $\overline\iota(X) \le 0$, thus $i = \iota(M_i) \le 0$,
thus $M_i$ is a sink module, say with center $p$. Note that $M_i$ is the quasi-top of $X$, 
thus $X' = M_i$. It is sufficient to show that $B(X) = B(M_i)$ (because this implies that $X$
is a sink module). Thus, we have to show that $T(F) \subseteq B(M_i)$ for all quasi-composition 
factors $F$ of $X$.

First, let us assume that we deal with an even component.
Since $M_i$ is a sink module, 
$i = -2t$ for some $t\ge 0$. Since $\overline\iota(X) \le 0$, 
the number of sink factors is greater than the
number of remaining factors, thus $i+2(l-1) \le 2t$. Thus, the quasi-length $l$ of $X$
is at most $2t+1$ and the quasi-composition factors of $X$ are 
$$
 M_{-2t}, M_{-2t+2},M_{-2t+4},\dots.
$$
Let $r$ be the radius of $M_0$. Then by Theorem 1, the radius of $M_{-2s}$ with $0\le s \le t$
is $r+2s\le r+2t$; in particular: the radius of $M_i$ is $r+2t$. Thus, for $0\le s \le t$,
$T(M_{-2s}) \subseteq B_{r+2t}(p).$ 
In general, if $T(M) \subseteq B_u(p)$, then $T(\tau M) \subseteq B_{u+2}(p)$. Thus, since
$T(M_0) \subseteq B_r(p)$ we see by induction that for $1\le s \le t$ we have
$T(M_{2s}) \subseteq B_{r+2s}(p) \subseteq B_{r+2t}(p) = B(M_i)$.
This shows that $T(F) \subseteq B(M_i)$ for all
factors $F$.

Second, assume that the component is odd.
Since $M_i$ is a sink module, 
$i = -2t-1$ for some $t\ge 0$. Since $\overline\iota(X) \le 0$, the number of sink factors is greater or equal the
number of remaining factors, therefore $i+2(l-1) \le 2t+1$. Thus, the quasi-length $l$ of $M$
is at most $2t+2$ and the quasi-composition factors of $M$ are 
$$
 M_{-2t-1}, M_{-2t+1}, M_{-2t+3},\dots.
$$
Again, let $r$ be the radius of $M_0$. Then by Theorem 1, the radius of $M_{-2s-1}$ 
with $0\le s \le t$
is $r+2s+1\le r+2t+1$; in particular: $B(M_i) = B_{r+2t+1}(p)$. 
Again using Theorem 1, we know that $T(M_1) \subseteq B_r(p)$. It follows that
$T(M_{2s+1}) \subseteq B_{r+2s}(p) \subseteq B_{r+2t+1}(p) = B(M_i)$.
This shows that $T(F) \subseteq B(M_i)$ for all factors $F$.
\hfill$\square$
	\bigskip 
We have found in this way many sink modules in a given regular component.
Here is a visualization:
$$
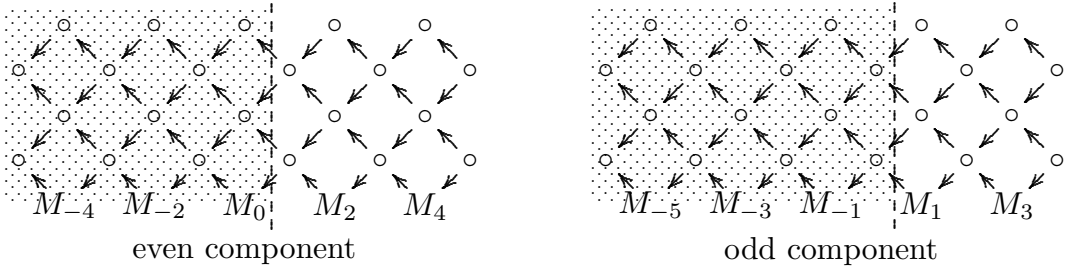
\hbox{\beginpicture
  \setcoordinatesystem units <.5cm,.5cm>
\put{\beginpicture
\put{$M_{-4}$} at 1 0
\put{$M_{-2}$} at 3 0 
\put{$M_{0}$} at 5 0 
\put{$M_{2}$} at 7 0 
\put{$M_{4}$} at 9 0 
\multiput{$\circ$} at 0 1  2 1  4 1  6 1  8 1  10 1
     1 2  3 2  5 2  7 2  9 2 
     0 3  2 3  4 3  6 3  8 3  10 3
     1 4  3 4  5 4  7 4  9 4 /
\arr{0.6 0.4}{0.3 0.7}
\arr{1.7 0.7}{1.4 0.4}
\arr{2.6 0.4}{2.3 0.7}
\arr{3.7 0.7}{3.4 0.4}
\arr{4.6 0.4}{4.3 0.7}
\arr{5.7 0.7}{5.4 0.4}
\arr{6.6 0.4}{6.3 0.7}
\arr{7.7 0.7}{7.4 0.4}
\arr{8.6 0.4}{8.3 0.7}
\arr{9.7 0.7}{9.4 0.4}

\arr{0.7 1.7}{0.3 1.3}
\arr{1.7 1.3}{1.3 1.7}
\arr{2.7 1.7}{2.3 1.3}
\arr{3.7 1.3}{3.3 1.7}
\arr{4.7 1.7}{4.3 1.3}
\arr{5.7 1.3}{5.3 1.7}
\arr{6.7 1.7}{6.3 1.3}
\arr{7.7 1.3}{7.3 1.7}
\arr{8.7 1.7}{8.3 1.3}
\arr{9.7 1.3}{9.3 1.7}

\arr{0.7 2.3}{0.3 2.7}
\arr{1.7 2.7}{1.3 2.3}
\arr{2.7 2.3}{2.3 2.7}
\arr{3.7 2.7}{3.3 2.3}
\arr{4.7 2.3}{4.3 2.7}
\arr{5.7 2.7}{5.3 2.3}
\arr{6.7 2.3}{6.3 2.7}
\arr{7.7 2.7}{7.3 2.3}
\arr{8.7 2.3}{8.3 2.7}
\arr{9.7 2.7}{9.3 2.3}

\arr{0.7 3.7}{0.3 3.3}
\arr{1.7 3.3}{1.3 3.7}
\arr{2.7 3.7}{2.3 3.3}
\arr{3.7 3.3}{3.3 3.7}
\arr{4.7 3.7}{4.3 3.3}
\arr{5.7 3.3}{5.3 3.7}
\arr{6.7 3.7}{6.3 3.3}
\arr{7.7 3.3}{7.3 3.7}
\arr{8.7 3.7}{8.3 3.3}
\arr{9.7 3.3}{9.3 3.7}
\setdashes <1mm>
\plot 5.6 -0.5  5.6 4.5 /
\put{even component} at 5 -1
\setshadegrid span <.7mm>
\vshade -.4 0 4.4  5.6 0 4.4 /
\endpicture} at 0 0
\put{\beginpicture
\put{$M_{-5}$} at 1 0
\put{$M_{-3}$} at 3 0 
\put{$M_{-1}$} at 5 0 
\put{$M_{1}$} at 7 0 
\put{$M_{3}$} at 9 0 
\multiput{$\circ$} at 0 1  2 1  4 1  6 1  8 1  10 1
     1 2  3 2  5 2  7 2  9 2 
     0 3  2 3  4 3  6 3  8 3  10 3
     1 4  3 4  5 4  7 4  9 4 /
\arr{0.6 0.4}{0.3 0.7}
\arr{1.7 0.7}{1.4 0.4}
\arr{2.6 0.4}{2.3 0.7}
\arr{3.7 0.7}{3.4 0.4}
\arr{4.6 0.4}{4.3 0.7}
\arr{5.7 0.7}{5.4 0.4}
\arr{6.6 0.4}{6.3 0.7}
\arr{7.7 0.7}{7.4 0.4}
\arr{8.6 0.4}{8.3 0.7}
\arr{9.7 0.7}{9.4 0.4}

\arr{0.7 1.7}{0.3 1.3}
\arr{1.7 1.3}{1.3 1.7}
\arr{2.7 1.7}{2.3 1.3}
\arr{3.7 1.3}{3.3 1.7}
\arr{4.7 1.7}{4.3 1.3}
\arr{5.7 1.3}{5.3 1.7}
\arr{6.7 1.7}{6.3 1.3}
\arr{7.7 1.3}{7.3 1.7}
\arr{8.7 1.7}{8.3 1.3}
\arr{9.7 1.3}{9.3 1.7}

\arr{0.7 2.3}{0.3 2.7}
\arr{1.7 2.7}{1.3 2.3}
\arr{2.7 2.3}{2.3 2.7}
\arr{3.7 2.7}{3.3 2.3}
\arr{4.7 2.3}{4.3 2.7}
\arr{5.7 2.7}{5.3 2.3}
\arr{6.7 2.3}{6.3 2.7}
\arr{7.7 2.7}{7.3 2.3}
\arr{8.7 2.3}{8.3 2.7}
\arr{9.7 2.7}{9.3 2.3}

\arr{0.7 3.7}{0.3 3.3}
\arr{1.7 3.3}{1.3 3.7}
\arr{2.7 3.7}{2.3 3.3}
\arr{3.7 3.3}{3.3 3.7}
\arr{4.7 3.7}{4.3 3.3}
\arr{5.7 3.3}{5.3 3.7}
\arr{6.7 3.7}{6.3 3.3}
\arr{7.7 3.3}{7.3 3.7}
\arr{8.7 3.7}{8.3 3.3}
\arr{9.7 3.3}{9.3 3.7}
\setdashes <1mm>
\plot 6.4 -0.5  6.4 4.5 /
\put{odd component} at 5 -1
\setshadegrid span <.7mm>
\vshade -.4 0 4.4  6.4 0 4.4 /
\endpicture} at 13 0
\endpicture}
$$
All the modules on the left of the dashes line are sink modules with center $p$.
	\bigskip 
Here is the dual assertion of Lemma 3. If $M$ is indecomposable and regular, with
a quasi-composition factor $F$. If the $\sigma$-orbit of $F$ contains $b$ flow modules, we
write $b(X) = b$ (note that this does not depend on the choice of $F$). 
	\medskip
{\bf Lemma 3*.} {\it Let $X$ be indecomposable and regular with $\overline\iota(X) \ge b(X)+1$.
Then $X$ is a source module and $B(X) = B(X'')$, where $X''$ is the quasi-socle of $X$.}
\hfill$\square$
	\medskip 
Here are the pictures dual to the previous ones; 
the modules on the right of the dashes line are source modules with a fixed
center, say $q$.
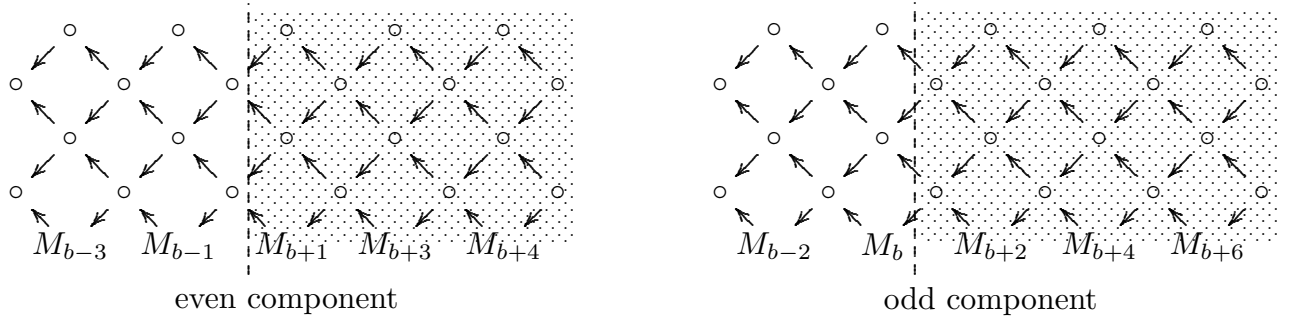
$$
\hbox{\beginpicture
  \setcoordinatesystem units <.6cm,.6cm>
\put{\beginpicture
\put{$M_{b-3}$} at 1 0
\put{$M_{b-1}$} at 3 0 
\put{$M_{b+1}$} at 5.1 0 
\put{$M_{b+3}$} at 7 0 
\put{$M_{b+4}$} at 9 0 
\multiput{$\circ$} at 0 1  2 1  4 1  6 1  8 1  10 1
     1 2  3 2  5 2  7 2  9 2 
     0 3  2 3  4 3  6 3  8 3  10 3
     1 4  3 4  5 4  7 4  9 4 /
\arr{0.6 0.4}{0.3 0.7}
\arr{1.7 0.7}{1.4 0.4}
\arr{2.6 0.4}{2.3 0.7}
\arr{3.7 0.7}{3.4 0.4}
\arr{4.6 0.4}{4.3 0.7}
\arr{5.7 0.7}{5.4 0.4}
\arr{6.6 0.4}{6.3 0.7}
\arr{7.7 0.7}{7.4 0.4}
\arr{8.6 0.4}{8.3 0.7}
\arr{9.7 0.7}{9.4 0.4}

\arr{0.7 1.7}{0.3 1.3}
\arr{1.7 1.3}{1.3 1.7}
\arr{2.7 1.7}{2.3 1.3}
\arr{3.7 1.3}{3.3 1.7}
\arr{4.7 1.7}{4.3 1.3}
\arr{5.7 1.3}{5.3 1.7}
\arr{6.7 1.7}{6.3 1.3}
\arr{7.7 1.3}{7.3 1.7}
\arr{8.7 1.7}{8.3 1.3}
\arr{9.7 1.3}{9.3 1.7}

\arr{0.7 2.3}{0.3 2.7}
\arr{1.7 2.7}{1.3 2.3}
\arr{2.7 2.3}{2.3 2.7}
\arr{3.7 2.7}{3.3 2.3}
\arr{4.7 2.3}{4.3 2.7}
\arr{5.7 2.7}{5.3 2.3}
\arr{6.7 2.3}{6.3 2.7}
\arr{7.7 2.7}{7.3 2.3}
\arr{8.7 2.3}{8.3 2.7}
\arr{9.7 2.7}{9.3 2.3}

\arr{0.7 3.7}{0.3 3.3}
\arr{1.7 3.3}{1.3 3.7}
\arr{2.7 3.7}{2.3 3.3}
\arr{3.7 3.3}{3.3 3.7}
\arr{4.7 3.7}{4.3 3.3}
\arr{5.7 3.3}{5.3 3.7}
\arr{6.7 3.7}{6.3 3.3}
\arr{7.7 3.3}{7.3 3.7}
\arr{8.7 3.7}{8.3 3.3}
\arr{9.7 3.3}{9.3 3.7}
\setdashes <1mm>
\plot 4.3 -0.5  4.3 4.5 /
\put{even component} at 5 -1
\setshadegrid span <.7mm>
\vshade 4.3 0 4.4  10.4 0 4.4 /

\endpicture} at 0 0
\put{\beginpicture
\put{$M_{b-2}$} at 1 0
\put{$M_{b}$} at 3 0 
\put{$M_{b+2}$} at 5 0 
\put{$M_{b+4}$} at 7 0 
\put{$M_{b+6}$} at 9 0 
\multiput{$\circ$} at 0 1  2 1  4 1  6 1  8 1  10 1
     1 2  3 2  5 2  7 2  9 2 
     0 3  2 3  4 3  6 3  8 3  10 3
     1 4  3 4  5 4  7 4  9 4 /
\arr{0.6 0.4}{0.3 0.7}
\arr{1.7 0.7}{1.4 0.4}
\arr{2.6 0.4}{2.3 0.7}
\arr{3.7 0.7}{3.4 0.4}
\arr{4.6 0.4}{4.3 0.7}
\arr{5.7 0.7}{5.4 0.4}
\arr{6.6 0.4}{6.3 0.7}
\arr{7.7 0.7}{7.4 0.4}
\arr{8.6 0.4}{8.3 0.7}
\arr{9.7 0.7}{9.4 0.4}

\arr{0.7 1.7}{0.3 1.3}
\arr{1.7 1.3}{1.3 1.7}
\arr{2.7 1.7}{2.3 1.3}
\arr{3.7 1.3}{3.3 1.7}
\arr{4.7 1.7}{4.3 1.3}
\arr{5.7 1.3}{5.3 1.7}
\arr{6.7 1.7}{6.3 1.3}
\arr{7.7 1.3}{7.3 1.7}
\arr{8.7 1.7}{8.3 1.3}
\arr{9.7 1.3}{9.3 1.7}

\arr{0.7 2.3}{0.3 2.7}
\arr{1.7 2.7}{1.3 2.3}
\arr{2.7 2.3}{2.3 2.7}
\arr{3.7 2.7}{3.3 2.3}
\arr{4.7 2.3}{4.3 2.7}
\arr{5.7 2.7}{5.3 2.3}
\arr{6.7 2.3}{6.3 2.7}
\arr{7.7 2.7}{7.3 2.3}
\arr{8.7 2.3}{8.3 2.7}
\arr{9.7 2.7}{9.3 2.3}

\arr{0.7 3.7}{0.3 3.3}
\arr{1.7 3.3}{1.3 3.7}
\arr{2.7 3.7}{2.3 3.3}
\arr{3.7 3.3}{3.3 3.7}
\arr{4.7 3.7}{4.3 3.3}
\arr{5.7 3.3}{5.3 3.7}
\arr{6.7 3.7}{6.3 3.3}
\arr{7.7 3.3}{7.3 3.7}
\arr{8.7 3.7}{8.3 3.3}
\arr{9.7 3.3}{9.3 3.7}
\setdashes <1mm>
\plot 3.6 -0.5  3.6 4.5 /
\put{odd component} at 5 -1
\setshadegrid span <.7mm>
\vshade 3.6 0 4.4  10.4 0 4.4 /

\endpicture} at 13 0
\endpicture}
$$

Now we apply Proposition 2 in order to see that the a module $X$ with
$1 \le \overline\iota(X) \le b(X)$
have to be flow modules. 
	\medskip
{\bf Proposition 4.} {\it For any indecomposable regular module $X$, we have $\overline\iota(X)
= \iota(X).$}
	\bigskip 
Proof. Let $\Cal X$ be a $\tau$-orbit. In order to show that $\overline\iota = \iota$ on
$\Cal X$, it is sufficient to to show that $\overline\iota(X) =
\iota(X)$ for some $X$ in $\Cal X$. Let $p$ be the center of the sink modules in $\Cal X$ and 
$q$ the center of the source modules in $\Cal X.$ Let $b = d(p,q)$.

First, assume that $b = 2s+1$ is odd. Then $\Cal X$ contains either $s$ or $s+1$ flow modules.
If $X$ belongs to $\Cal X$ and $\overline\iota(X) = 0,$ then $X$ is a sink module, and
the modules $\tau X,\dots,\tau^sX$ are the flow modules in $\Cal X$. Thus $\Cal X$
has precisely $s$ flow modules and according to Proposition 2, $X$ has index $0$. 
Thus $\overline\iota(X) = 0 = \iota X$. 

If $X$ belongs to $\Cal X$ and $\overline\iota(X) = -1,$ then $X$ is a sink module, and
the modules $\tau X,\dots,\tau^{s+1}X$ are the flow modules in $\Cal X$. Thus $\Cal X$
has precisely $s+1$ flow modules and again using Proposition 2, we see that $X$ has index $-1$. 
Thus $\overline\iota(X) = -1 = \iota X$. 

Second, assume that $b = 2s$ is even, thus $\Cal X$ contains precisely $s$ flow modules.
We denote by $M_i$ quasi-composition factors of the modules in $\Cal X$, with
$\iota(M_i) = i$.  
Let $l$ be the quasi-length of the modules in $\Cal X$. Assume that $X$ belongs to $\Cal X$
with $\overline\iota(X) = 0.$ Then $X$ is a sink module, 
the modules $\tau X,\dots,\tau^sX$ are the flow modules in $\Cal X$, and $\tau^{s+1} X$
is a source module. The quasi-top of $X$ is $M_{-l}$, 
the quasi-socle of $\tau^{s+1}X$ is $M_{b+1+l+1}$. According to Lemma 3, 
$B(X) = B(M_{-l})$, thus $r(M) = r_0+l$
According to Lemma 3*, $B(\tau M^{s+1}) = B(M_{b+1+l+1})$, thus $r(\tau M^{s+1}) = r_0+l+1$.
This shows that $r(M) < r(\tau M^{s+1})$. Thus, Proposition 3 asserts that $\iota(X) = 0.$
A similar calculation shows that for $\overline\iota(X) = -1,$ we have $\iota(X) = -1.$

Thus, in all cases we have found a module $X\in \Cal X$ with $\overline\iota(X) = \iota(X)$
and therefore $\overline\iota(X') = \iota(X')$ for all $X'\in \Cal X$.
  \hfill$\square$
	\medskip 
The distribution of the sink modules, flow modules and source modules in the regular
Auslander-Reiten components is as follows:
$$
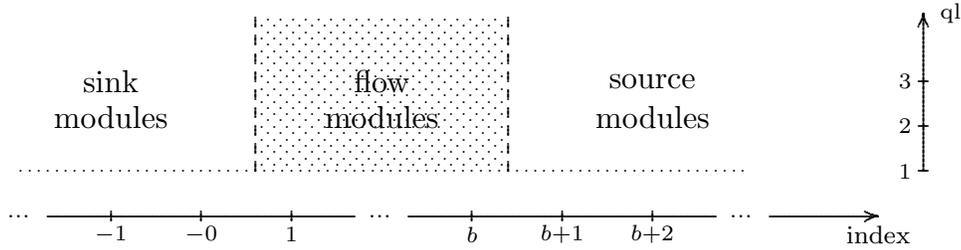
\hbox{\beginpicture
  \setcoordinatesystem units <1cm,.5cm>
\put{} at 0 4
\put{sink} at  -3 2
\put{modules} at  -3 1.2
\put{flow} at   0 2
\put{modules} at  0 1.2
\put{source} at  3 2
\put{modules} at  3 1.2
\setdots <1mm>
\plot -4 0  4 0 /
\setdashes <1mm>
\plot -1.4 0  -1.4 3.5 /
\plot 1.4 0  1.4 3.5 /
\put{$\ssize -1$} at  -3 -1.4
\put{$\ssize -0$} at  -2 -1.4  
\put{$\ssize 1$} at  -1 -1.4
\multiput{$\ssize \dots$} at 0 -1 -4 -1  4 -1 /  
  
\put{$\ssize b$} at  1 -1.4
\put{$\ssize b+1$} at  2 -1.4
\put{$\ssize b+2$} at  3 -1.4
\setshadegrid span <.7mm>
\vshade -1.4 0 3.4  1.4 0 3.4 /

\setsolid
\arr{4.3 -1}{5.5 -1}
\plot -3.7 -1 -.3 -1 /
\plot  0.3 -1 3.7 -1 /
\put{\rmk index} at 5.5 -1.4 
\plot -3 -1.1 -3 -.9 /
\plot -2 -1.1 -2 -.9 /
\plot -1 -1.1 -1 -.9 /
\plot  1 -1.1 1 -.9 /
\plot  2 -1.1 2 -.9 /
\plot  3 -1.1 3 -.9 /

\arr{6 0}{6 3.5}
\plot 5.95 0  6.05 0 /
\plot 5.95 1  6.05 1 /
\plot 5.95 2  6.05 2 /
\put{$\ssize 1$} at 5.8 0
\put{$\ssize 2$} at 5.8 1
\put{$\ssize 3$} at 5.8 2
\put{$\ssize \ql$} at 6.3 3.5
\endpicture}
$$
Of course, Proposition 4 implies (and actually is equivalent to) the assertion of
Theorem 3. 
	\medskip
Let us look now at Theorem 4. Actually, it remains to deal with the last assertion of Theorem 4.
Clearly, it is sufficient to show this for
one module in the $\tau$-orbit of $Y$, thus we may assume that $\iota(Y) \le 0$.
Let $M$ be the quasi-top of $Y$. According to Lemma 3, we have $B(Y) = B(M)$, 
thus $r(Y) = r(M)$.
Recall that $\eta Y = Y',$ where
$0\to X \to Y\oplus Y' \to Z \to 0$ is an Auslander-Reiten sequence with $|Y| < |Y'|.$
Then $\tau^-M$ is the quasi-top both of $Y'$ and $Z$. Therefore
$B(\eta Y) = B(Y') = B(\tau^-M)$ and we know that $r(\tau^-M) = r(M)+2,$ thus $r(\eta Y) = 
r(Y') = r(M)+2 = r(Y)+2.$
Concerning the centers, the equalities $B(\eta Y) = B(\tau^-M)$ and $B(Y) = B(M)$ yield
$C(\eta Y) = C(\tau^-M) = C(M) = C(Y)$. 
	\hfill$\square$
	\bigskip\bigskip
We have seen: If $M,M'$ belong to the same Auslander-Reiten component,
then
$$
 p(M) = p(M'),\quad q(M) = q(M'),\quad r_0(M)-\ql(M) = r_0(M')-\ql(M').
$$
Thus, if $\Cal C$ is a regular component, and $M$ belongs to $\Cal C$, we may define
$$
\align
 p(\Cal C) &= p(M),\cr
 q(\Cal C) &= q(M),\cr
 b(\Cal C) &= d(p(M),q(M)),\cr
 r(\Cal C) &= r_0(M)-\ql(M).
\endalign
$$
The radius $r(M)$ of a module $M$ in $\Cal C$ is 
$$
 r(\Cal C) + R_{b(\Cal C)}(\iota(M),\ql (M)), 
$$
where $R_b:\Bbb Z\times \Bbb N_1 \to \Bbb N_0$ 
is the following function
$$
R_b(i,l) = \left\{\matrix -i+l & &i\le 0,\quad\ \ \cr
               l-1 &\text{for}& 1\le i\le b,\cr
               i+l& &b<i.\quad\  \ 
             \endmatrix \right. 
$$
At the end of the paper, we exhibit the function $R_b(i,l)$ for some special values of $b$.
	\bigskip\bigskip
{\bf 8. A warning, two questions and several remarks.}
     \medskip
{\bf The Warning.} Let $M$ be an indecomposable regular module. Its
center path $\pi(M)$ may not be contained in the support $T(M)$ of
$M$ (but it is always contained in the ball $B(M)$). Here is a
typical example: Let $\pi = (q = a_0, a_1,\dots, a_b = q)$ be a path
with $q$ a source. Then $P_0(p)$ is a submodule of $P_b(q)$
and $q$ is not contained in the support of $M = P_b(p)/P_0(q)$.
This is the module $M = M_1$ constructed as Case I in the proof of
Theorem 2, see section 4. There, we have seen that $\pi(M)$ is
the given path $\pi.$
	\bigskip
{\bf Remark 1.} The case $n= 2$.
We have seen that for $n\ge 3$ (thus for the wild Kronecker algebras) 
almost all modules in the $\tau$-orbit of a regular graded module are 
sink or source modules. In contrast, in the tame case $n=2$, the sink modules are the
preprojective modules, the source modules are the preinjective modules, thus all
the indecomposable regular graded modules are flow modules. We should stress that
also for $n=2$, the regular components are of the form
$\Bbb Z\Bbb A_\infty$ (and not proper quotients), so that any regular component
contains infinitely many quasi-simple modules. To repeat: The number of
quasi-simple flow modules in a regular component is finite, if $n\ge 3$, but it is infinite,
if $n = 2.$ 
As in our investigation of bristles see Appendix C of [R4], 
we encouter a finiteness result for wild cases 
which is not valid in the corresponding tame case.
	\medskip 

{\bf Remark 2.} When dealing with wild hereditary algebras, one knows that
the regular modules have exponential growth when applying $\tau$ (or $\tau^-$),
see [R1,B,K1,K2]. This concerns the dimension vector (and all its
coefficients). In contrast, the main observation of this paper asserts that
the radius of an indecomposable regular module eventually grows linearly when 
we iterate the application of $\tau$: 
it increases step by step by 2. Of course, we have seen also exponential growth (at the end
of section 4): the number of boundary vertices in the support of a sink module increases 
exponentially, when we apply $\tau$. Note that the number $\beta(M)$ of boundary vertices
is a lower bound for the length of $M$.
	\medskip 
{\bf Remark 3.} Our main motivation for these investigations is our interest
in (ungraded) Kronecker modules: to find invariants for regular Kronecker modules.
The present paper provides such invariants for the gradable modules (those which are
obtained from graded Kronecker modules by forgetting the grading, or, equivalently,
for those representations of $K(n)$ which can be lifted to its universal covering
$(T(n),\Omega)$). Given an indecomposable 
gradable Kronecker module $\overline M$ with cover module
$M$, the numerical invariants $r_0(M)$ and $b(M)$ are invariants of $\overline M$.
Actually, it is $r_0$ and 
the path $\pi(M)$ (or better, its equivalence class under the covering
group) which are the decisive invariants of $\overline M$.
 	\medskip 
{\bf Question 1.} We have stressed in remark 2 that dealing with graded Kronecker
modules, the radius of an indecomposable regular module eventually increases 
linearly when we iterate the application of $\tau$. 
Is there a similar invariant for ungraded Kronecker modules with this growth behaviour?
	\medskip 
{\bf Remark 4.} If $\Cal C$ is a regular component in $\mo(T(n),\Omega)$,
iterated applications of $\tau$ to any module in $\Cal C$
provide balls with a fixed center $p$ (in particular, $p$ is an invariant of $\Cal C$).
The preprojective and the preinjective 
component behave completely different, there is no focus on a special
vertex of $T(n).$ 
	\medskip 
{\bf Remark 5.} In this paper, when looking at a representations $M$ 
of $(T(n),\Omega)$, we actually were dealing just
with the dimension vector of $M$, not with $M$ itself. Indecomposable representations of
$(T(n),\Omega)$ with same dimension vectors have the same behavior with respect to
the action of $\sigma$: If $M,M'$ are indecomposable representations of $(T(n),\Omega)$
with $\bdim M = \bdim M'$,
then $p(M) = p(M'), q(M) = q(M')$ and $r_0(M) = r_0(M').$
	\medskip 
{\bf Remark 6.} Why is the index of a module defined in such a way that it 
decreases along irreducible maps?
In general, we like to draw abelian categories by focusing the attention to
the direction of maps, drawing arrows from left to right. 
According to Baer [B] and Kerner [K1,K2], there is a  
global direction of a hereditary 
module category, defined by the maps from the projective modules to the injective
modules (going from left to right): for wild hereditary algebras, this global
direction is opposite to the direction of the irreducible maps (thus, opposite to $\tau$).
To repeat, dealing with regular modules of wild hereditary algebras, the irreducible maps 
point in the opposite direction of the global direction of the category. Thus, it is reasonable
to draw the regular Auslander-Reiten components for a hereditary
algebra by using arrows from right to left.
	\medskip 
{\bf Question 2.}
The definition of the index of an indecomposable representation of 
$(T(n),\Omega)$ relies not only on $\mo(T(n),\Omega)$, but also on 
$\mo(T(n),\sigma\Omega)$. Thus, let $M$ be an indecomposable regular $(T(n),\Omega)$-module 
say with index $i$. Then $\tau M$ is an indecomposable regular $(T(n),\Omega)$-module 
say with index $i-2$. But attached to $M$ are two modules with index $i-1$, namely
the  $(T(n),\Omega)$-module $\mu(M)$ 
which arises as the middle term of the Auslander-Reiten sequence
ending in $M$ (this module is not necessarily indecomposable, it may be the direct
sum of two indecomposable modules with index $i-1$), and the $(T(n),\sigma\Omega)$-module 
$\sigma M.$ What is the precise relationship between $\mu(M)$ and $\sigma M$? 
	\medskip
{\bf Remark 7.} The topic of this paper may be compared to the so-called
Game of Life, as introduced by Conway, see [Gar], as a cellular automaton.
Whereas the
Game of Life is played on the grid with set of vertices $\Bbb Z^2$,
here we start with the $n$-regular tree $T(n)$. If we want to increase
the analogy, we should deal not with dimension vectors of modules
but work modulo $2$, thus dealing with functions $f\:T(n)_0 \to \{0,1\}$
with finite support.

	\bigskip\bigskip
{\bf 9. References.}
     \medskip
\item{[B]} D.~Baer: Wild hereditary Artin algebras and linear methods.
   manuscripta mathematica 55, 69–-82.

\item{[BGP]} I.~N.~Bernstein, I.~M.~Gelfand, V.~A.~Ponomarev:
      Coxeter functors and Gabriel's theorem. Uspekhi Mat. Nauk 28 (1973),
      Russian Math. Surveys 28 (1973), 17--32.

\item{[BB]} S.~Brenner, M.~C.~R.~Butler: The equivalence of certain functors occuring
 in the representation theory of artin algebras and species. J\. London Math\. Soc\.
  14 (1976), 183--187.

\item{[FR1]} Ph.~Fahr, C.~M.~Ringel: A partition formula for Fibonacci numbers.
    Journal of Integer Sequences, Vol. 11 (2008), Article 08.1.4.

\item{[FR2]} Ph.~Fahr, C.~M.~Ringel: Categorification of the Fibonacci Numbers 
   Using Representations of Quivers
    Journal of Integer Sequences. Vol. 15 (2012), Article 12.2.1.

\item{[Gab]} P.~Gabriel: Auslander-Reiten sequences and representation-finite algebras. In:
 Representation Theory I.
  Lecture Notes in Math., 831, Berlin, New York: Springer-Verlag, (1980), 1--71.

\item{[Gar]} M.~Gardner: Martin: Mathematical Games. The fantastic combinations of
John Conway's new solitaire game {\it life}. Scientific American 223 (October 1970), 120--123.

\item{[K1]} O.~Kerner: Representations of wild quivers. In: Representation Theory 
   of Algebras and Related Topics. CMS Conference Proceedings 19. Amer. Math.
   Soc. (1994), 65-107
\item{[K2]} O.~Kerner: More representations of wild quivers. Contemp. Math. 607 (2014), 35-55

\item{[R1]} C.~M.~Ringel: Representations of K-species and bimodules. 
  J. Algebra 41 (1976), 269-302. 

\item{[R2]} C.~M.~Ringel: Finite-dimensional hereditary algebras of wild representation type. 
  Math. Z. 161 (1978), 235-255. 

\item{[R3]} C.~M.~Ringel: Infinite dimensional representations of finite dimensional hereditary 
algebras. Symposia Math. 23 (1979), 321-412. 

\item{[R4]}  C.~M.~Ringel: Kronecker modules generated by modules of length 2.
    To appear in: Representations of Algebras. 
    Amer.Math.Soc. Contemporary Mathematics. \newline
   arXiv:1612.07679. 

\bigskip\bigskip
{\rmk Claus Michael Ringel

Fakult\"at f\"ur Mathematik, Universit\"at Bielefeld

D-33501 Bielefeld, Germany\par}

{\ttk ringel\@math.uni-bielefeld.de}
\vfill\eject
{\bf The function $R_b(i,l)$ on even and an odd components, for some values $b$.}
	\medskip 
$$
\hbox{\beginpicture
  \setcoordinatesystem units <.4cm,.39cm>
\put{\beginpicture
\put{$b=0$} [l] at -6.5 8.5
\put{\rmk index $\ssize i$\strut} at 9 7
\setdots <1mm>
\multiput{$1$} at -1 0 /
\multiput{$2$} at  0 1  1 0 /
\multiput{$3$} at -3 0  -2 1  -1 2  /
\multiput{$4$} at 0 3  1 2  2 1  3 0 /
\multiput{$5$} at  -5 0  -4 1  -3 2  -2 3  -1 4 /
\multiput{$6$} at  0 5  1 4  2 3  3 2  4 1  5 0 /
\multiput{$7$} at  -6 1  -5 2  -4 3  -3 4  -2 5  /
\multiput{$8$} at 2 5  3 4  4 3  5 2  6 1  /
\multiput{$9$} at  -6 3  -5 4  -4 5 / 
\multiput{$10$} at 4 5  5 4  6 3 /
\multiput{$11$} at  -6 5 /
\multiput{$12$} at  6 5 /

\plot -6 1  -5 0  0 5  5 0  6 1 /
\plot -6 1  -2 5  3 0   6 3   /
\plot -6 3  -4 5  1 0  6 5 /
\plot -6 3  -3 0  2 5  6 1 /
\plot -6 5  -1 0  4 5  6 3 /

\plot -6.8 .2  -6 1  -6.8 1.8 /
\plot -6.8 2.2  -6 3  -6.8 3.8 /
\plot -6.8 4.2  -6 5  -6.8 5.8 /

\plot  -6 5  -5.2 5.8 /
\plot -4.8 5.8  -4 5  -3.2 5.8 /
\plot -2.8 5.8  -2 5  -1.2 5.8 /
\plot  6 5  5.2 5.8 /
\plot 4.8 5.8  4 5  3.2 5.8 /
\plot 2.8 5.8  2 5  1.2 5.8 /
\plot 0.8 5.8  0 5  -.8 5.8 /

\plot 6.8 0.2  6 1  6.8 1.8 /
\plot 6.8 2.2  6 3  6.8 3.8 /
\plot 6.8 4.2  6 5  6.8 5.8 /

\setdashes <1mm>
\plot -0.5 -.5  -.5 5.8 /
\plot  -.6 -.5   -.6 5.8 /

\put{$\ssize -5$\strut} at -6 7
\put{$\ssize -4$\strut} at -5 7
\put{$\ssize -3$\strut} at -4 7
\put{$\ssize -2$\strut} at -3 7
\put{$\ssize -1$\strut} at -2 7
\put{$\ssize 0$\strut} at -1 7
\put{$\ssize 1$\strut} at 0 7
\put{$\ssize 2$\strut} at 1 7
\put{$\ssize 3$\strut} at 2 7
\put{$\ssize 4$\strut} at 3 7
\put{$\ssize 5$\strut} at 4 7
\put{$\ssize 6$\strut} at 5 7
\put{$\ssize 7$\strut} at 6 7
\put{\rmk even component} at 0 -1.5
\setdots <.6mm>
\plot -6.7 -.1  -5.3 -.1 /
\plot -4.7 -.1  -3.3 -.1 /
\plot -2.7 -.1  -1.3 -.1 /
\plot -.7 -.1  0.7 -.1 /
\plot 6.7 -.1  5.3 -.1 /
\plot 4.7 -.1  3.3 -.1 /
\plot 2.7 -.1  1.3 -.1 /
\endpicture} at 0 0
\put{\beginpicture
\put{\strut} at -5.5 8.5
\setdots <1mm>
\multiput{$1$} at 0 0 /
\multiput{$2$} at -2 0  -1 1 /
\multiput{$3$} at  0 2  1 1  2 0 /
\multiput{$4$} at -4 0  -3 1  -2 2  -1 3 / 
\multiput{$5$} at 0 4  1 3  2 2  3 1  4 0 /
\multiput{$6$} at -6 0  -5 1  -4 2  -3 3  -2 4  -1 5 /
\multiput{$7$} at  1 5  2 4  3 3  4 2  5 1  6 0 /
\multiput{$8$} at -6 2  -5 3  -4 4  -3 5 /
\multiput{$9$} at  3 5  4 4  5 3  6 2 /
\multiput{$10$} at -6 4  -5 5 /
\multiput{$11$} at  5 5  6 4 /
\plot       -6 0  -1 5  4 0  6 2 /
\plot -6 2  -4 0  1 5  6 0   /
\plot -6 4  -2 0  3 5  6 2 /
\plot -6 2  -3 5  2 0  6 4 /
\plot -6 4  -5 5  0 0  5 5  6 4 /

\plot -6.8 0.8  -6 0 /
\plot -6.8 1.2  -6 2  -6.8 2.8 /
\plot -6.8 3.2  -6 4  -6.8 4.8 /

\plot -5.8 5.8  -5 5  -4.2 5.8 /
\plot -3.8 5.8  -3 5  -2.2 5.8 /
\plot -1.8 5.8  -1 5  -0.2 5.8 /
\plot 5.8 5.8  5 5  4.2 5.8 /
\plot 3.8 5.8  3 5  2.2 5.8 /
\plot 1.8 5.8  1 5  0.2 5.8 /

\plot 6.8 0.8  6 0 /
\plot 6.8 1.2  6 2  6.8 2.8 /
\plot 6.8 3.2  6 4  6.8 4.8 /

\setdashes <1mm>
\plot -0.5 -.5  -.5 5.8 /
\plot -.6 -.5   -.6 5.8 /
\put{$\ssize -5$\strut} at -6 7
\put{$\ssize -4$\strut} at -5 7
\put{$\ssize -3$\strut} at -4 7
\put{$\ssize -2$\strut} at -3 7
\put{$\ssize -1$\strut} at -2 7
\put{$\ssize 0$\strut} at -1 7
\put{$\ssize 1$\strut} at 0 7
\put{$\ssize 2$\strut} at 1 7
\put{$\ssize 3$\strut} at 2 7
\put{$\ssize 4$\strut} at 3 7
\put{$\ssize 5$\strut} at 4 7
\put{$\ssize 6$\strut} at 5 7
\put{$\ssize 7$\strut} at 6 7
\put{\rmk odd component} at 0 -1.5
\setdots <.6mm>
\plot -6.8 -.1  -6.3 -.1 /
\plot  6.8 -.1   6.3 -.1 /
\plot -5.7 -.1  -4.3 -.1 /
\plot -3.7 -.1  -2.3 -.1 /
\plot -1.7 -.1  -0.3 -.1 /
\plot 5.7 -.1  4.3 -.1 /
\plot 3.7 -.1  2.3 -.1 /
\plot 1.7 -.1  0.3 -.1 /

\endpicture} at 16 0

\endpicture}
$$

$$
\hbox{\beginpicture
  \setcoordinatesystem units <.4cm,.39cm>
\put{\beginpicture
\put{$b=1$} [l] at -6.5 8.5
\put{\rmk index $\ssize i$\strut} at 9 7
\setdots <1mm>
\multiput{$1$} at -1 0  0 1  1 0 /
\multiput{$3$} at -3 0  -2 1  -1 2  0 3  1 2  2 1  3 0 /
\multiput{$5$} at  -5 0  -4 1  -3 2  -2 3  -1 4  0 5  1 4  2 3  3 2  4 1  5 0 /
\multiput{$7$} at  -6 1  -5 2  -4 3  -3 4  -2 5  2 5  3 4  4 3  5 2  6 1  /
\multiput{$9$} at  -6 3  -5 4  -4 5  4 5  5 4  6 3 /
\multiput{$11$} at  -6 5  6 5 /

\plot -6 1  -5 0  0 5  5 0  6 1 /
\plot -6 1  -2 5  3 0   6 3   /
\plot -6 3  -4 5  1 0  6 5 /
\plot -6 3  -3 0  2 5  6 1 /
\plot -6 5  -1 0  4 5  6 3 /

\plot -6.8 .2  -6 1  -6.8 1.8 /
\plot -6.8 2.2  -6 3  -6.8 3.8 /
\plot -6.8 4.2  -6 5  -6.8 5.8 /

\plot  -6 5  -5.2 5.8 /
\plot -4.8 5.8  -4 5  -3.2 5.8 /
\plot -2.8 5.8  -2 5  -1.2 5.8 /
\plot  6 5  5.2 5.8 /
\plot 4.8 5.8  4 5  3.2 5.8 /
\plot 2.8 5.8  2 5  1.2 5.8 /
\plot 0.8 5.8  0 5  -.8 5.8 /

\plot 6.8 0.2  6 1  6.8 1.8 /
\plot 6.8 2.2  6 3  6.8 3.8 /
\plot 6.8 4.2  6 5  6.8 5.8 /
\setdashes <1mm>
\plot -0.5 -.5  -.5 5.8 /
\plot   .5 -.5   .5 5.8 /
\setshadegrid span <.7mm>
\vshade -.5 -.5 5.8 .5 -.5 5.8 /
\put{$\ssize -5$\strut} at -6 7
\put{$\ssize -4$\strut} at -5 7
\put{$\ssize -3$\strut} at -4 7
\put{$\ssize -2$\strut} at -3 7
\put{$\ssize -1$\strut} at -2 7
\put{$\ssize 0$\strut} at -1 7
\put{$\ssize 1$\strut} at 0 7
\put{$\ssize 2$\strut} at 1 7
\put{$\ssize 3$\strut} at 2 7
\put{$\ssize 4$\strut} at 3 7
\put{$\ssize 5$\strut} at 4 7
\put{$\ssize 6$\strut} at 5 7
\put{$\ssize 7$\strut} at 6 7
\put{\rmk even component} at 0 -1.5
\setdots <.6mm>
\plot -6.7 -.1  -5.3 -.1 /
\plot -4.7 -.1  -3.3 -.1 /
\plot -2.7 -.1  -1.3 -.1 /
\plot -.7 -.1  0.7 -.1 /
\plot 6.7 -.1  5.3 -.1 /
\plot 4.7 -.1  3.3 -.1 /
\plot 2.7 -.1  1.3 -.1 /
\endpicture} at 0 0

\put{\beginpicture
\put{\strut} at -5.5 8.5
\setdots <1mm>
\multiput{$0$} at 0 0 /
\multiput{$2$} at -2 0  -1 1  0 2  1 1  2 0 /
\multiput{$4$} at -4 0  -3 1  -2 2  -1 3  0 4  1 3  2 2  3 1  4 0 /
\multiput{$6$} at -6 0  -5 1  -4 2  -3 3  -2 4  -1 5  1 5  2 4  3 3  4 2  5 1  6 0 /
\multiput{$8$} at -6 2  -5 3  -4 4  -3 5  3 5  4 4  5 3  6 2 /
\multiput{$10$} at -6 4  -5 5  5 5  6 4 /
\plot       -6 0  -1 5  4 0  6 2 /
\plot -6 2  -4 0  1 5  6 0   /
\plot -6 4  -2 0  3 5  6 2 /
\plot -6 2  -3 5  2 0  6 4 /
\plot -6 4  -5 5  0 0  5 5  6 4 /

\plot -6.8 0.8  -6 0 /
\plot -6.8 1.2  -6 2  -6.8 2.8 /
\plot -6.8 3.2  -6 4  -6.8 4.8 /

\plot -5.8 5.8  -5 5  -4.2 5.8 /
\plot -3.8 5.8  -3 5  -2.2 5.8 /
\plot -1.8 5.8  -1 5  -0.2 5.8 /
\plot 5.8 5.8  5 5  4.2 5.8 /
\plot 3.8 5.8  3 5  2.2 5.8 /
\plot 1.8 5.8  1 5  0.2 5.8 /

\plot 6.8 0.8  6 0 /
\plot 6.8 1.2  6 2  6.8 2.8 /
\plot 6.8 3.2  6 4  6.8 4.8 /
\setdashes <1mm>
\plot -0.5 -.5  -.5 5.8 /
\plot   .5 -.5   .5 5.8 /
\setshadegrid span <.7mm>
\vshade -.5 -.5 5.8 .5 -.5 5.8 /
\put{$\ssize -5$\strut} at -6 7
\put{$\ssize -4$\strut} at -5 7
\put{$\ssize -3$\strut} at -4 7
\put{$\ssize -2$\strut} at -3 7
\put{$\ssize -1$\strut} at -2 7
\put{$\ssize 0$\strut} at -1 7
\put{$\ssize 1$\strut} at 0 7
\put{$\ssize 2$\strut} at 1 7
\put{$\ssize 3$\strut} at 2 7
\put{$\ssize 4$\strut} at 3 7
\put{$\ssize 5$\strut} at 4 7
\put{$\ssize 6$\strut} at 5 7
\put{$\ssize 7$\strut} at 6 7
\put{\rmk odd component} at 0 -1.5
\setdots <.6mm>
\plot -6.8 -.1  -6.3 -.1 /
\plot  6.8 -.1   6.3 -.1 /
\plot -5.7 -.1  -4.3 -.1 /
\plot -3.7 -.1  -2.3 -.1 /
\plot -1.7 -.1  -0.3 -.1 /
\plot 5.7 -.1  4.3 -.1 /
\plot 3.7 -.1  2.3 -.1 /
\plot 1.7 -.1  0.3 -.1 /

\endpicture} at 16 0

\endpicture}
$$

$$
\hbox{\beginpicture
  \setcoordinatesystem units <.4cm,.39cm>
\put{\beginpicture
\put{$b=4$} [l] at -6.5 8.5
\put{\rmk index $\ssize i$\strut} at 9 7

\setdots <1mm>
\multiput{$0$} at 1 0  3 0 /
\multiput{$1$} at -1 0  0 1  2 1 /
\multiput{$2$} at 1 2  3 2   4 1  5 0 /
\multiput{$3$} at -3 0  -2 1  -1 2  0 3  2 3  /
\multiput{$4$} at 1 4     /
\multiput{$5$} at  -5 0  -4 1  -3 2  -2 3  -1 4  0 5  2 5 /
\multiput{$4$} at   3 4  4 3  5 2  6 1  /
\multiput{$7$} at  -6 1  -5 2  -4 3  -3 4  -2 5  /
\multiput{$6$} at 4 5  5 4  6 3  /
\multiput{$9$} at  -6 3  -5 4  -4 5 / 
\multiput{$8$} at 6 5   /
\multiput{$11$} at  -6 5 /

\plot -6 1  -5 0  0 5  5 0  6 1 /
\plot -6 1  -2 5  3 0   6 3   /
\plot -6 3  -4 5  1 0  6 5 /
\plot -6 3  -3 0  2 5  6 1 /
\plot -6 5  -1 0  4 5  6 3 /

\plot -6.8 .2  -6 1  -6.8 1.8 /
\plot -6.8 2.2  -6 3  -6.8 3.8 /
\plot -6.8 4.2  -6 5  -6.8 5.8 /

\plot  -6 5  -5.2 5.8 /
\plot -4.8 5.8  -4 5  -3.2 5.8 /
\plot -2.8 5.8  -2 5  -1.2 5.8 /
\plot  6 5  5.2 5.8 /
\plot 4.8 5.8  4 5  3.2 5.8 /
\plot 2.8 5.8  2 5  1.2 5.8 /
\plot 0.8 5.8  0 5  -.8 5.8 /

\plot 6.8 0.2  6 1  6.8 1.8 /
\plot 6.8 2.2  6 3  6.8 3.8 /
\plot 6.8 4.2  6 5  6.8 5.8 /
\setdashes <1mm>
\plot -0.5 -.5  -.5 5.8 /
\plot  3.5 -.5   3.5 5.8 /
\setshadegrid span <.7mm>
\vshade -.5 -.5 5.8  3.5 -.5 5.8 /
\put{$\ssize -5$\strut} at -6 7
\put{$\ssize -4$\strut} at -5 7
\put{$\ssize -3$\strut} at -4 7
\put{$\ssize -2$\strut} at -3 7
\put{$\ssize -1$\strut} at -2 7
\put{$\ssize 0$\strut} at -1 7
\put{$\ssize 1$\strut} at 0 7
\put{$\ssize 2$\strut} at 1 7
\put{$\ssize 3$\strut} at 2 7
\put{$\ssize 4$\strut} at 3 7
\put{$\ssize 5$\strut} at 4 7
\put{$\ssize 6$\strut} at 5 7
\put{$\ssize 7$\strut} at 6 7
\put{\rmk even component} at 0 -1.5
\setdots <.6mm>
\plot -6.7 -.1  -5.3 -.1 /
\plot -4.7 -.1  -3.3 -.1 /
\plot -2.7 -.1  -1.3 -.1 /
\plot -.7 -.1  0.7 -.1 /
\plot 6.7 -.1  5.3 -.1 /
\plot 4.7 -.1  3.3 -.1 /
\plot 2.7 -.1  1.3 -.1 /

\endpicture} at 0 0
\put{\beginpicture
\put{\strut} at -5.5 8.5
\setdots <1mm>
\multiput{$0$} at 0 0  2 0 /
\multiput{$2$} at -2 0  -1 1  0 2  2 2 /
\multiput{$1$} at  1 1  3 1  4 0  /
\multiput{$4$} at -4 0  -3 1  -2 2  -1 3  0 4  2 4 / 
\multiput{$3$} at  1 3  3 3  4 2  5 1  6 0  /
\multiput{$6$} at -6 0  -5 1  -4 2  -3 3  -2 4  -1 5 /
\multiput{$5$} at  1 5  3 5   4 4  5 3  6 2  /
\multiput{$8$} at -6 2  -5 3  -4 4  -3 5 /
\multiput{$10$} at -6 4  -5 5 /
\multiput{$7$} at  5 5  6 4 /
\plot       -6 0  -1 5  4 0  6 2 /
\plot -6 2  -4 0  1 5  6 0   /
\plot -6 4  -2 0  3 5  6 2 /
\plot -6 2  -3 5  2 0  6 4 /
\plot -6 4  -5 5  0 0  5 5  6 4 /

\plot -6.8 0.8  -6 0 /
\plot -6.8 1.2  -6 2  -6.8 2.8 /
\plot -6.8 3.2  -6 4  -6.8 4.8 /

\plot -5.8 5.8  -5 5  -4.2 5.8 /
\plot -3.8 5.8  -3 5  -2.2 5.8 /
\plot -1.8 5.8  -1 5  -0.2 5.8 /
\plot 5.8 5.8  5 5  4.2 5.8 /
\plot 3.8 5.8  3 5  2.2 5.8 /
\plot 1.8 5.8  1 5  0.2 5.8 /

\plot 6.8 0.8  6 0 /
\plot 6.8 1.2  6 2  6.8 2.8 /
\plot 6.8 3.2  6 4  6.8 4.8 /
\setdashes <1mm>
\plot -0.5 -.5  -.5 5.8 /
\plot 3.5 -.5   3.5 5.8 /
\setshadegrid span <.7mm>
\vshade -.5 -.5 5.8  3.5 -.5 5.8 /
\put{$\ssize -5$\strut} at -6 7
\put{$\ssize -4$\strut} at -5 7
\put{$\ssize -3$\strut} at -4 7
\put{$\ssize -2$\strut} at -3 7
\put{$\ssize -1$\strut} at -2 7
\put{$\ssize 0$\strut} at -1 7
\put{$\ssize 1$\strut} at 0 7
\put{$\ssize 2$\strut} at 1 7
\put{$\ssize 3$\strut} at 2 7
\put{$\ssize 4$\strut} at 3 7
\put{$\ssize 5$\strut} at 4 7
\put{$\ssize 6$\strut} at 5 7
\put{$\ssize 7$\strut} at 6 7
\put{\rmk odd component} at 0 -1.5
\setdots <.6mm>
\plot -6.8 -.1  -6.3 -.1 /
\plot  6.8 -.1   6.3 -.1 /
\plot -5.7 -.1  -4.3 -.1 /
\plot -3.7 -.1  -2.3 -.1 /
\plot -1.7 -.1  -0.3 -.1 /
\plot 5.7 -.1  4.3 -.1 /
\plot 3.7 -.1  2.3 -.1 /
\plot 1.7 -.1  0.3 -.1 /

\endpicture} at 16 0

\endpicture}
$$

$$
\hbox{\beginpicture
  \setcoordinatesystem units <.4cm,.39cm>
\put{\beginpicture
\put{$b=5$} [l] at -6.5 8.5
\put{\rmk index $\ssize i$\strut} at 9 7
\setdots <1mm>
\multiput{$0$} at   1 0  3 0 /
\multiput{$1$} at -1 0  0 1    2 1  4 1  5 0 /
\multiput{$2$} at   1 2  3 2 /
\multiput{$3$} at -3 0  -2 1  -1 2  0 3   2 3  4 3  5 2  6 1 /
\multiput{$4$} at   1 4  3 4  /
\multiput{$5$} at 0 5   2 5  /
\multiput{$5$} at   -5 0  -4 1  -3 2  -2 3  -1 4     4 5  5 4  6 3    /
\multiput{$7$} at  -6 1  -5 2  -4 3  -3 4  -2 5  6 5 /

\multiput{$9$} at   -6 3  -5 4  -4 5  /
\multiput{$11$} at  -6 5  /

\plot -6 1  -5 0  0 5  5 0  6 1 /
\plot -6 1  -2 5  3 0   6 3   /
\plot -6 3  -4 5  1 0  6 5 /
\plot -6 3  -3 0  2 5  6 1 /
\plot -6 5  -1 0  4 5  6 3 /

\plot -6.8 .2  -6 1  -6.8 1.8 /
\plot -6.8 2.2  -6 3  -6.8 3.8 /
\plot -6.8 4.2  -6 5  -6.8 5.8 /

\plot  -6 5  -5.2 5.8 /
\plot -4.8 5.8  -4 5  -3.2 5.8 /
\plot -2.8 5.8  -2 5  -1.2 5.8 /
\plot  6 5  5.2 5.8 /
\plot 4.8 5.8  4 5  3.2 5.8 /
\plot 2.8 5.8  2 5  1.2 5.8 /
\plot 0.8 5.8  0 5  -.8 5.8 /

\plot 6.8 0.2  6 1  6.8 1.8 /
\plot 6.8 2.2  6 3  6.8 3.8 /
\plot 6.8 4.2  6 5  6.8 5.8 /
\setdashes <1mm>
\plot -0.5 -.5  -.5 5.8 /
\plot   4.5 -.5   4.5 5.8 /
\setshadegrid span <.7mm>
\vshade -.5 -.5 5.8  4.5 -.5 5.8 /
\put{$\ssize -5$\strut} at -6 7
\put{$\ssize -4$\strut} at -5 7
\put{$\ssize -3$\strut} at -4 7
\put{$\ssize -2$\strut} at -3 7
\put{$\ssize -1$\strut} at -2 7
\put{$\ssize 0$\strut} at -1 7
\put{$\ssize 1$\strut} at 0 7
\put{$\ssize 2$\strut} at 1 7
\put{$\ssize 3$\strut} at 2 7
\put{$\ssize 4$\strut} at 3 7
\put{$\ssize 5$\strut} at 4 7
\put{$\ssize 6$\strut} at 5 7
\put{$\ssize 7$\strut} at 6 7
\put{\rmk even component} at 0 -1.5
\setdots <.6mm>
\plot -6.7 -.1  -5.3 -.1 /
\plot -4.7 -.1  -3.3 -.1 /
\plot -2.7 -.1  -1.3 -.1 /
\plot -.7 -.1  0.7 -.1 /
\plot 6.7 -.1  5.3 -.1 /
\plot 4.7 -.1  3.3 -.1 /
\plot 2.7 -.1  1.3 -.1 /

\endpicture} at 0 0

\put{\beginpicture
\put{\strut} at -5.5 8.5
\setdots <1mm>
\multiput{$0$} at 0 0  2 0  4 0 /
\multiput{$1$} at 1 1  3 1 /
\multiput{$3$} at 1 3  3 3 /
\multiput{$5$} at 1 5  3 5 /
\multiput{$2$} at -2 0  -1 1  0 2  2 2   4 2  5 1  6 0   /
\multiput{$4$} at -4 0  -3 1  -2 2  -1 3  0 4  2 4   4 4  5 3  6 2 /
\multiput{$6$} at -6 0  -5 1  -4 2  -3 3  -2 4  -1 5  /
\multiput{$8$} at -6 2  -5 3  -4 4  -3 5  5 5  6 4 /
\multiput{$10$} at -6 4  -5 5 /
\plot       -6 0  -1 5  4 0  6 2 /
\plot -6 2  -4 0  1 5  6 0   /
\plot -6 4  -2 0  3 5  6 2 /
\plot -6 2  -3 5  2 0  6 4 /
\plot -6 4  -5 5  0 0  5 5  6 4 /

\plot -6.8 0.8  -6 0 /
\plot -6.8 1.2  -6 2  -6.8 2.8 /
\plot -6.8 3.2  -6 4  -6.8 4.8 /

\plot -5.8 5.8  -5 5  -4.2 5.8 /
\plot -3.8 5.8  -3 5  -2.2 5.8 /
\plot -1.8 5.8  -1 5  -0.2 5.8 /
\plot 5.8 5.8  5 5  4.2 5.8 /
\plot 3.8 5.8  3 5  2.2 5.8 /
\plot 1.8 5.8  1 5  0.2 5.8 /

\plot 6.8 0.8  6 0 /
\plot 6.8 1.2  6 2  6.8 2.8 /
\plot 6.8 3.2  6 4  6.8 4.8 /
\setdashes <1mm>
\plot -0.5 -.5  -.5 5.8 /
\plot   4.5 -.5   4.5 5.8 /
\setshadegrid span <.7mm>
\vshade -.5 -.5 5.8  4.5 -.5 5.8 /
\put{$\ssize -5$\strut} at -6 7
\put{$\ssize -4$\strut} at -5 7
\put{$\ssize -3$\strut} at -4 7
\put{$\ssize -2$\strut} at -3 7
\put{$\ssize -1$\strut} at -2 7
\put{$\ssize 0$\strut} at -1 7
\put{$\ssize 1$\strut} at 0 7
\put{$\ssize 2$\strut} at 1 7
\put{$\ssize 3$\strut} at 2 7
\put{$\ssize 4$\strut} at 3 7
\put{$\ssize 5$\strut} at 4 7
\put{$\ssize 6$\strut} at 5 7
\put{$\ssize 7$\strut} at 6 7
\put{\rmk odd component} at 0 -1.5
\setdots <.6mm>
\plot -6.8 -.1  -6.3 -.1 /
\plot  6.8 -.1   6.3 -.1 /
\plot -5.7 -.1  -4.3 -.1 /
\plot -3.7 -.1  -2.3 -.1 /
\plot -1.7 -.1  -0.3 -.1 /
\plot 5.7 -.1  4.3 -.1 /
\plot 3.7 -.1  2.3 -.1 /
\plot 1.7 -.1  0.3 -.1 /

\endpicture} at 16 0

\endpicture}
$$

\bye